\documentclass[ijoc,nonblindrev]{./StyleFiles/informs3arXiv}

\OneAndAHalfSpacedXI 

\usepackage{./StyleFiles/endnotes}
\usepackage[skip=15pt]{subcaption}

\usepackage[noend]{algpseudocode}
\usepackage{etoolbox}

\makeatletter
\newcommand*{\algrule}[1][\algorithmicindent]{%
  \makebox[#1][l]{%
    \hspace*{.7em}
    \color{gray}\vrule height .75\baselineskip depth .25\baselineskip
  }
}

\newcount\ALG@printindent@tempcnta
\def\ALG@printindent{%
    \ifnum \theALG@nested>0
    \ifx\ALG@text\ALG@x@notext
    \else
    \unskip
    \ALG@printindent@tempcnta=1
    \loop
    \algrule[\csname ALG@ind@\the\ALG@printindent@tempcnta\endcsname]%
    \advance \ALG@printindent@tempcnta 1
    \ifnum \ALG@printindent@tempcnta<\numexpr\theALG@nested+1\relax
    \repeat
    \fi
    \fi
}
\patchcmd{\ALG@doentity}{\noindent\hskip\ALG@tlm}{\ALG@printindent}{}{\errmessage{failed to patch}}
\patchcmd{\ALG@doentity}{\item[]\nointerlineskip}{}{}{} 
\makeatother

\newcounter{subroutine}
\makeatletter
\newenvironment{subroutine}[1][htb]{%
  \let\c@algorithm\c@subroutine
  \renewcommand{\ALG@name}{Subroutine}
  \begin{algorithm}[#1]%
  }{\end{algorithm}
}
\makeatother

\let\footnote=\endnote
\let\enotesize=\normalsize
\def\notesname{Endnotes}%
\def\makeenmark{$^{\theenmark}$}
\def\enoteformat{\rightskip0pt\leftskip0pt\parindent=1.75em
  \leavevmode\llap{\theenmark.\enskip}}

\usepackage{./StyleFiles/natbib}
 \bibpunct[, ]{(}{)}{,}{a}{}{,}%
 \def\bibfont{\small}%
 \def\bibsep{\smallskipamount}%
 \def\bibhang{24pt}%
\TheoremsNumberedThrough     
\ECRepeatTheorems
\EquationsNumberedThrough    

\usepackage{float}

\usepackage{graphicx}
\usepackage{./StyleFiles/custom_formats}
\usepackage{./StyleFiles/custom_ian}
\usetikzlibrary{patterns}
\hypersetup{
colorlinks,
linkcolor={red!50!black},
urlcolor={blue!80!black}
}

\newlength\myindent
\algdef{SE}[SUBALG]{Indent}{EndIndent}{}{\algorithmicend\ }%
\algtext*{Indent}
\algtext*{EndIndent}

\newcommand{\iz}[1]{\textcolor{black}{#1}}

\theoremstyle{plain}

\usepackage{longtable}
\usepackage{pdflscape}

\usepackage{tikz}
\usetikzlibrary{decorations.pathreplacing,calc}
\usetikzlibrary{shapes.geometric}

\usepackage[mathscr]{eucal}

\usepackage{soul,xcolor}

\begin{document}

\setstcolor{blue}

\RUNAUTHOR{Bodur, Chan and Zhu}

\RUNTITLE{Inverse Mixed Integer Linear Optimization}

\TITLE{Inverse Mixed Integer Optimization: \\ Polyhedral Insights and Trust Region Methods}


\ARTICLEAUTHORS{%
\AUTHOR{Merve Bodur, Timothy C. Y. Chan, Ian Yihang Zhu}
\AFF{Department of Mechanical and Industrial Engineering, University of Toronto, Toronto, Ontario M5S 3G8, Canada} 
} 

\ABSTRACT{Inverse optimization -- determining parameters of an optimization problem that render a given solution optimal -- has received increasing attention in recent years. While significant inverse optimization literature exists for convex optimization problems, there have been few advances  for discrete problems, despite the ubiquity of applications that fundamentally rely on discrete decision-making. In this paper, we present a new set of theoretical insights and algorithms for the general class of inverse mixed integer linear optimization problems. Specifically, a general characterization of optimality conditions is established and leveraged to design new cutting plane solution algorithms. Through an extensive set of computational experiments, we show that our methods provide substantial improvements over existing methods in solving the largest and most difficult instances to date.}


\maketitle



\section{Introduction}
Inverse optimization -- inferring unobserved parameters of a (\emph{forward}) optimization problem that render a given (\emph{forward-feasible}) solution optimal -- has received increasing attention in recent years. Most of the literature has focused on inverse optimization for convex forward optimization problems, allowing the bi-level inverse problem to be reformulated as a single-level convex problem using the Karush-Kuhn-Tucker conditions. Examples include linear \citep[e.g.,][]{zhang1996calculating, ahuja2001inverse}, multiobjective linear \citep[e.g.][]{naghavi2019inverse}, conic \citep[e.g.][]{iyengar2005inverse} and general convex optimization problems \citep[e.g.][]{zhang2010inverse, zhang2010perturbation}. 
\iz{These models have been applied across a wide range of application domains (see Section \ref{subsubsec:applications}) and have also been extended methodologically} 
in several modern directions
including estimating model parameters from multiple observed solutions while considering concepts such as statistical consistency \citep{aswani2018inverse} or distributional robustness \citep{esfahani2018data}. 

In contrast, there have been very few advances in inverse optimization for discrete \iz{forward} problems, despite the ubiquity of methodological research in discrete optimization and real-world applications that have discrete decisions. Part of the challenge may stem from the fact that the optimality conditions for discrete problems like mixed integer linear optimization (MILO) generally do not lead to computationally tractable solution algorithms like they do for convex problems. For example, \cite{schaefer2009inverse} and \cite{lamperski2015polyhedral} demonstrate how superadditive duality can be used to develop exact reformulations of bi-level inverse integer and mixed-integer optimization problems, respectively. However, the resulting single-level formulation is an exponentially large linear program that is intractable beyond problems with few variables. 

The other main idea in the literature for a general-purpose inverse MILO solution method is a cutting plane algorithm. 
First proposed in \cite{wang2009cutting}, the idea is to decompose the inverse MILO problem into a master problem and a subproblem. The former is a relaxation of the inverse MILO problem that provides candidates for parameters to be inferred, whereas the latter identifies \emph{extreme points} of the convex hull of the forward-feasible region that generate cuts to send back to the master problem. While it was shown that small problems could be solved efficiently, this approach does not scale well to larger problems. The computational cost of computing each extreme point is large in general, since it involves solving an instance of the forward MILO problem. The only advance to this method since its development was a proposed heuristic that parallelizes the computation of extreme points \citep{duan2011heuristic}.

In this paper, we develop a new cutting plane framework for solving inverse MILO problems. 
We begin by providing a more general characterization of the optimality conditions of the inverse MILO problem relative to previous literature.
We demonstrate that extreme points are sufficient but not necessary for characterizing inverse optimality, and that more parsimonious representations of the problem exist. Using this insight, we develop a new cutting plane method where the cut generation subroutine efficiently identifies interior forward-feasible points using trust regions. These cuts come at significantly lower computational cost compared to cuts generated from extreme points. In addition, we observe that inverse MILO problems can be solved with many fewer cuts compared to the classical cutting plane algorithm. 

Our specific contributions are: 

\begin{enumerate}
    \item We present a novel characterization of optimality conditions for inverse MILO problems by introducing the concept of \emph{generator sets} (Section \ref{sec:backgroundandmotivation}). Generator sets characterize the complete family of reformulations of the MILO problem that preserve the feasible region of the inverse optimization problem. Our main result provides necessary and sufficient conditions for any set to be a generator set. We also show that there exist ``small'' generator sets, which are attractive from a computational perspective.

    \item Motivated by our theoretical insights, we propose a new family of cutting plane methods that restrict the feasible region of the MILO problem using trust regions (Section \ref{sec:general_framework}). We also propose computational enhancements that further speed up the algorithm (Section \ref{sec:enhancements}). Our algorithms generalize the classical cutting plane approach to solving inverse MILO problems.
    
    \item We propose an extension of the inverse MILO problem to a setting with multiple input data points and illustrate how our cutting plane algorithm can be easily extended to solve this problem (Section \ref{sec:contemporary_models}). 
    
    \item Through a comprehensive numerical study \iz{over instances drawn from the MIPLIB 2017 benchmark library \citep{gleixner2019miplib}}, we demonstrate that our new cutting plane algorithm significantly outperforms the state-of-the-art algorithm by solving the largest and most difficult problem instances in the literature to date (Sections \ref{sec:numerical_setup} and \ref{sec:numerical_ex}). We observe that our new cut generation methods can both substantially speed up cut generation and reduce the number of cuts required. 

\end{enumerate}

\subsection{Related literature}\label{subsec:literature_review}

Below, we review relevant inverse optimization models and solution methods, as well as the main application areas in which these models \iz{can be} found.

\subsubsection{Models and methods.}
The classic\iz{al} inverse optimization problem, which is to infer a set of model parameters that render a given feasible solution optimal, has been studied over a wide range of problem settings. The parameters that are inferred can be the cost vector or the constraints. While there are a few studies that focus on estimating constraint parameters \citep{guler2010capacity, birge2017inverse, chan2020inverse}, the vast majority of papers focus on estimating the cost vector. \iz{Most early works focused on deriving algorithms to solve inverse models for particular forward problems (e.g., shortest path, minimum spanning tree), while more recent works consider general-purpose methods for broader classes of forward problems (e.g., mixed integer optimization). \cite{heuberger2004inverse} provides an overview of the early works whereas the literature discussed in the introduction are examples of more recent works.} \iz{In this paper, we focus on estimating the cost vector for general mixed integer forward problems, particularly in the form of the classical inverse model}.

Inverse optimization problems with ``noisy" data have also been studied \citep[e.g.][]{troutt2006behavioral, keshavarz2011imputing, chan2014generalized, chan2019inverse}. These problems are characterized by the fact that there does not exist a set of (non-trivial) parameters that render a given solution, or a set of solutions, optimal. For these problems, the cost vector is estimated by minimizing various notions of suboptimality \citep[e.g.][]{bertsimas2015data, aswani2018inverse, esfahani2018data, babier2021ensemble}. \iz{Although we focus predominantly on the classical inverse model in the paper,}  we show that the insights and methods developed can be easily extended to a model for multiple ``noisy" data points that resembles the models found for inverse convex optimization.

Finally, we borrow the term ``trust region" from the general optimization community, a term which broadly describes the restriction of a particular search space. The use of trust region concepts have appeared in a range of different domains such as nonlinear optimization and, more relevant to our paper, decomposition algorithms. Within the latter domain, trust regions have been applied to the master problems of decomposition models for stochastic programming, and are shown to help generate stronger cuts in this context \citep{linderoth2003decomposition, zverovich2012computational, rahmaniani2017benders}. In contrast, we apply trust regions to the cut generation problem in our decomposition framework, which leads to improvements in both strength of cuts \emph{and} cut generation time over a wide array of problem structures.

\subsubsection{Applications.}\label{subsubsec:applications}

Inverse optimization has been studied in a wide variety of applications \iz{to estimate latent parameters and infer subjective preferences using observed decision data. 
Examples of such applications can be found in} energy markets \citep{ruiz2013revealing, saez2016data, birge2017inverse}, healthcare \citep{erkin2010eliciting, chan2014generalized,babier2020knowledge}, finance \citep{bertsimas2012inverse, utz2014tri, yu2020learning}, and transportation \citep{chow2012inverse, chow2014nonlinear, xu2018network, zhang2018price}, \iz{where inverse optimization is used to provide insights into market structure, treatment design, risk aversion and route-choice preferences, respectively. This literature makes use of both classical inverse models and inverse models for ``noisy" decision data.} 

\iz{Inverse optimization is also prevalent in the bi-level optimization literature, especially in regards to pricing and incentive design problems.}
\iz{Recent} examples of such applications include the design of tolls \iz{in traffic and commodity transportation networks} \citep{marcotte2009toll,brotcorne2011exact, esfandeh2016regulating, kuiteing2017network, kuiteing2018pricing}, \iz{price schedules and carbon taxes} in energy systems \citep{zhou2011designing, afcsar2021revenue}, and profit-sharing mechanisms in carrier alliances \citep{agarwal2010network, houghtalen2011designing}. 
\iz{Here, the classical inverse optimization model appears as a ``subproblem" within the solution process used to solve these pricing problems. Specifically, a ``master problem" generates a decision ``target" (e.g., a sustainable routing decision) while an inverse model computes a set of objective perturbations (e.g., tolls on roads) for which the target becomes an optimal decision for an underlying decision-maker (e.g., a transportation company).} \iz{The master problem can also propose new targets if the inverse subproblem is infeasible \citep{afcsar2021revenue}.} 

Despite the ubiquity of \iz{decision-making} models in these application areas that fundamentally rely on \iz{making discrete choices}, existing literature has focused primarily on inverse optimization for continuous forward problems. The few exceptions to date are in energy planning \citep{zhou2011designing}, transportation \citep{chow2012inverse}, and sustainability \citep{turner2013examining}. The first two applications were formulated as general inverse MILO problems and applied the classical cutting plane algorithm \citep{wang2009cutting}. The sustainability application, an inverse knapsack problem, was reformulated using superadditive duality and solved by restricting focus to linear functions. 

\subsection{Notation}

Throughout the paper, vectors and matrices are written in bold, while sets are defined using calligraphic letters. Subscripts denote specific elements of a vector, whereas superscripts denote different vectors. We use $(\cdot)^\top$ to denote the transpose operator. 
For any set $\mathcal{S}$, $\conv(\mathcal{S})$ is its convex hull, $\ext(\mathcal{S})$ is its set of extreme points, and $\mE(\mathcal{S}) = \ext(\conv(\mathcal{S}))$. 

\section{Structure of Inverse MILO Problems}\label{sec:backgroundandmotivation}

In this section, we study the structure of inverse MILO problems. We provide a novel characterization of inverse-feasibility over general non-convex forward optimization problems. We do this through the definition of a \emph{generator set}, which characterizes the complete family of reformulations of the forward optimization problem that preserve the feasible region of the inverse optimization problem. 
In turn, we show how certain feasible region definitions are amenable to more efficient solution methods for the inverse optimization problem. Proofs of our results are provided in the Appendix. 

\subsection{Problem definition}

Our inverse optimization model is based on the following 
MILO problem, known as the \emph{forward problem}: 
%
%
\begin{align}\label{model:FP}
\begin{split}
\textbf{FP}(\bc, \mX): \ \underset{\bx}{\text{minimize}} \quad & \bc^\top \bx \\
\text{subject to} \quad & \bx \in \mX := \{\bA \bx \geq \bb, \ \bx \in \mathbb{Z}^{n-q} \times \mathbb{R}^{q} \}.
\end{split}
\end{align}
%
Let $\mF(\bc,\mX)$ be the optimal solution set of $\textbf{FP}(\bc, \mX)$. Elements of $\mX$ and $\mF(\bc,\mX)$ are called \emph{forward-feasible} and \emph{forward-optimal}, respectively.
\iz{For a given forward-feasible solution $\bhx \in \mX$, let $\mC(\bhx,\mX)$ define the \emph{inverse-feasible region}, which is the set of all cost vectors $\bc \in \mathbb{R}^n$ that render the solution $\bhx$ optimal, i.e., 
\begin{equation}\label{eq:inverse-feasibility}
\mC(\bhx, \mX) := \{\bc \in \mathbb{R}^n \; | \; \bhx \in \mF(\bc, \mX)\}.
\end{equation}}
%
\iz{The inverse optimization problem is defined as the following optimization model,} 
\begin{align}
\label{model:IOM-Bilevel1a}
\begin{split}
\ \textbf{IO}(\bc^0,\hat\bx, \mX) : \ \underset{\bc \iz{\in \mP}}{\text{minimize}} \quad & g(\bc) \\ 
\text{subject to} \quad & \bc \in \mC(\hat\bx, \mX),
\end{split}
\end{align}
\iz{which seeks to find a cost vector $\bc$ that minimizes an objective function $g(\bc)$, is inverse-feasible, and satisfies any context-specific constraints $\bc \in \mP$. A common choice of $g(\bc)$ is%
\begin{align*}
g(\bc) := \norm{\bc-\bc^0}_1.
\end{align*}}
\iz{The inverse optimization model with this objective seeks to find feasible cost vectors $\bc$ that minimize the deviation to a given reference cost vector $\bc^0$, as measured by the Manhattan distance.} \iz{Practically, this objective function is relevant in many applications (see Section \ref{subsec:literature_review} for details). First, in estimation problems, $\bc^0$ can represent an initial estimate of the cost vector, and given an observed decision $\bhx$, the cost estimates can be refined (minimally) to be compatible with $\bhx$. For example, $\bc^0$ may represent an estimate of network link costs used to describe route preferences (e.g., where initial estimates are based on road distances), and an observed routing decision $\bhx$ can be used to refine these estimates \citep{burton1992instance, ahuja2001inverse, chen2021inverse}. Second, in bi-level optimization problems, $\bc^0$ represents the known cost vector of a decision-maker, and a set of objective perturbations to $\bc^0$ can be designed and prescribed such that a decision target $\bhx$ becomes an optimal choice for the decision-maker. For example, the elements of $\bc^0$ may represent the cost of producing energy using different resources, in which case the vector $\bc - \bc^0$ would describe the value of carbon taxes or subsidies that must be imposed on each resource for a low-emission production plan $\bhx$ to become optimal for an energy producer \citep{zhou2011designing, rathore2021differential}.} \iz{In these settings, there may be context-specific constraints on $\bc$, such as} non-negativity ($\mP = \mathbb{R}_+^{n}$), integrality ($\mP = \mathbb{Z}^n$), \iz{variable bounds ($\mP = [\mathbf{l}, \mathbf{u}]^n$)}, or a transformation into a multi-objective problem ($\mP = \{ \sum_{i = 1}^m \alpha_i \bv^i,  \ \alpha_i \geq 0, i = 1, \ldots, m\}$ where $\bv^1, \ldots, \bv^m$ are given objectives and $\boldsymbol{\alpha}$ is a decision vector of objective weights). 

\iz{From this point on, we present our inverse optimization model as
\begin{align}
\label{model:IOM-Bilevel1}
\begin{split}
\ \textbf{IO}(\bc^0,\hat\bx, \mX) : \ \underset{\bc \iz{\in \mP}}{\text{minimize}} \quad & \norm{\bc - \bc^0}_1 \\ 
\text{subject to} \quad & \bc \in \mC(\hat\bx, \mX).
\end{split}
\end{align}}
\iz{We make this modeling choice because the model is well-motivated and because it is the formulation that is studied in all previous works on inverse mixed integer optimization \citep{wang2009cutting, schaefer2009inverse, duan2011heuristic, lamperski2015polyhedral, bulutcomplexity}. Nonetheless, we emphasize that this modeling choice comes without loss of generality. The theoretical insights and methodological contributions presented in this paper focus exclusively on difficulties posed by the inverse-feasibility constraint (equation \eqref{eq:inverse-feasibility}), and are thus directly applicable for any inverse MILO problem in the form of model \eqref{model:IOM-Bilevel1a}}. \iz{We also assume without loss of generality that the inverse optimization problem is feasible. If it is not, then the methods proposed in this paper can prove its infeasibility.}

In the remainder of this section (Sections \ref{subsec:prelim} and \ref{subsec:opti_conditions}), we focus on the core difficulty of solving the inverse MILO problem: characterizing the inverse-feasible region $\mC(\bhx, \mX)$. In contrast to convex optimization problems, where this condition can be transformed into convex constraints using the Karush-Kuhn-Tucker conditions, determining whether a given $\hat\bx$ is optimal with respect to some $\bc$, let alone the $\bc$ that optimizes the objective function of the inverse problem, requires solving the forward MILO problem, which is NP-hard in general. \iz{Furthermore, attempting to use inverse linear optimization models to approximate inverse MILO models by ignoring the integrality constraints of the forward-feasible region can result in arbitrarily bad approximations. This is further discussed and illustrated in Section \ref{EC:approximations} of the Electronic Companion.}

\subsection{Preliminary remarks}\label{subsec:prelim}

We first note three basic properties of the inverse-feasible region, then describe it using a set of linear constraints.

\begin{remark}
The inverse-feasible region is nonempty, since $\textbf{0} \in \mC(\hat\bx, \mX)$.
\end{remark}

\begin{remark}
\iz{The inverse-feasible region contains only the zero vector, i.e., $\mC(\hat\bx, \mX) = \{\mathbf{0}\}$, if and only if $\mX$ is full dimensional and $\bhx$ is not on the boundary of conv($\mX$). If there exists at least one equality constraint or one binary variable used in the description of $\mX$, then $\mC(\bx, \mX) \neq \{\mathbf{0}\}$ for all $\bx \in \mX$. Similarly, if there exists at least one constraint that is binding at $\bhx$, including variable bounds, then $\mC(\hat\bx, \mX)  \neq \{\mathbf{0}\}$.}
\end{remark}

\begin{remark}
\label{rmk:Cbdd}
While the forward-feasible region $\mX$ does not have to be bounded, the inverse-feasible region $\mC(\hat\bx, \mX)$ will only contain objective vectors for which the forward problem is bounded. In particular, if $\bc^{rec}$ is a nonzero vector in the recession cone of $\conv(\mX)$, then we have $\bc^{rec} \notin \mC$. The reason is that $\ix$ is assumed to be forward-optimal with a finite objective value (under an unobserved objective function)  and $\inf \{ (\bc^{rec})^\top \bx \ | \ \bx \in \mX \} = -\infty$. Thus $\ix \notin \arg \min \{ (\bc^{rec})^\top \bx \ | \ \bx \in \mX \}$. In other words, no such $\bc^{rec}$ can make $\ix$ forward-optimal. 
\end{remark}

A natural way to explicitly formulate model~\eqref{model:IOM-Bilevel1} is to re-write $\bc \in \mC(\hat\bx,\mX)$ using a potentially infinite set of linear constraints %
\begin{equation}\label{cons:full_set}
\bc^\top(\hat\bx - \bx) \leq 0, \ \forall \bx \in \mX,
\end{equation}
suggesting that a cutting plane method can be used to solve $\mathbf{IO}(\bc^0, \hat\bx,\mX)$. Indeed, \cite{wang2009cutting} proposed a simple cutting plane algorithm to solve the inverse MILO problem in a finite (but exponential, in general) number of iterations by noting that it is sufficient to replace $\mX$ in equation \eqref{cons:full_set} with the extreme points of its convex hull. 

\begin{remark}{\citep{wang2009cutting}}\label{rmk:extreme_point_formulation}
Model~\eqref{model:IOM-Bilevel1} is equivalent to
\begin{subequations}
\label{model:IOM-Bilevel2}
\begin{align}
\ 
\underset{\bc \iz{\in \mP}}{\text{minimize}} \quad & \norm{\bc - \bc^0}_1 \label{obj:IOM-Bilevel2} \\
\text{subject to} \quad & \bc^\top(\hat\bx - \bx^j) \leq 0, \quad \forall \bx^j \in \mE(\mX). \label{con:IOM-Bilevel2} 
\end{align}
\end{subequations}
\end{remark}

The cutting plane algorithm proposed in \citet{wang2009cutting} starts without constraints \eqref{con:IOM-Bilevel2}, and iteratively adds cuts of this form through the computation of new extreme points $\bx^j$. The efficiency of this approach depends heavily on the difficulty of generating cuts, which amounts to solving the forward MILO problem, as well as the number of extreme points of the convex hull of $\mX$ \citep{bulutcomplexity}. With the exception of parallelizing the search for extreme points \citep{duan2011heuristic}, the algorithm proposed in \citet{wang2009cutting} remains the state-of-the-art in solving inverse MILO problems.

\subsection{Optimality conditions}\label{subsec:opti_conditions}

The main insight from \cite{wang2009cutting} is that replacing constraints \eqref{cons:full_set} with \eqref{con:IOM-Bilevel2} preserves the set of optimal solutions to the inverse optimization problem, i.e., $\mC(\bhx, \mE(\mX)) = \mC(\bhx, \mX)$. In this subsection, we characterize the full family of sets $\mG$ such that $\mC(\bhx,\mG) = \mC(\bhx,\mX)$. 

\begin{definition}
Given $\mX$, any set $\mG \subseteq \mathbb{R}^n$ for which 
\begin{equation}
\mC(\hat\bx, \mG) = \mC(\hat\bx, \mX)
\end{equation} 
is a \emph{generator set} for $\mX$. If $\mG \subseteq \mX$, then $\mG$ is a \emph{forward-feasible generator set}.
\end{definition}
Examples of generator sets include $\conv(\mX)$ and $\mE(\mX)$, with the latter being a forward-feasible generator set. 

Next, we demonstrate there may exist a large family of such sets for $\mX$. First, note that constraint set \eqref{cons:full_set} defines a polyhedral cone. Second, each half-space $\bc^\top(\hat\bx-\bx) \le 0$ \emph{induced by} a feasible point $\bx \in \mX$ is not uniquely \emph{determined by} $\bhx-\bx$, but rather by any multiple of the vector $\bhx- \bx$. In other words, the half-space induced by $\bx$ is the same as the one induced by $\hat\bx + \lambda (\bx - \hat\bx)$ for any $\lambda > 0$. The same cut can thus be induced by infinitely many other points of $\mathbb{R}^n$.

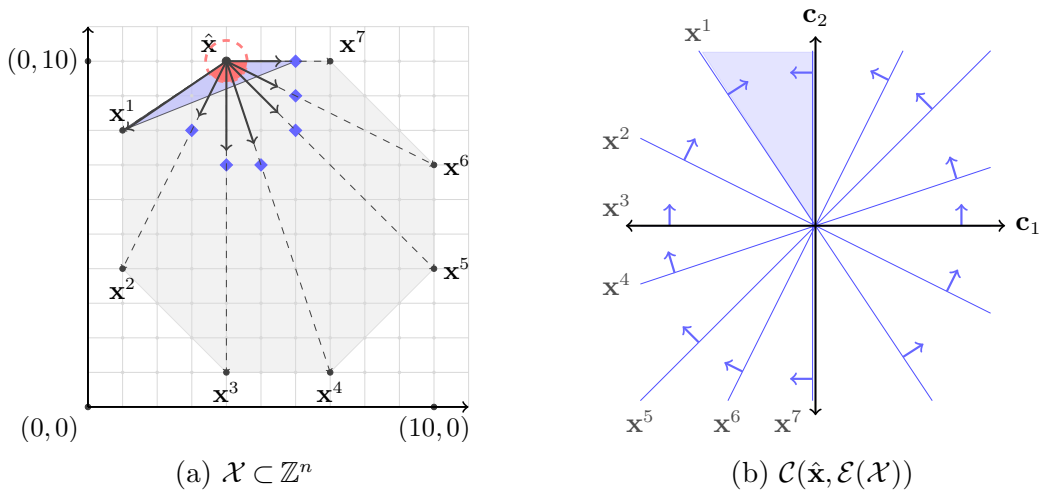
\begin{figure}[h]
\centering
\begin{subfigure}[t]{0.48\textwidth}
\centering
\scalebox{0.92}{
  \begin{tikzpicture}
    \begin{scope}
    \clip (-1.7,-1) rectangle (5.2cm,5.2cm); 
    \begin{scope}[transparency group]
    \begin{scope}[blend mode=multiply]
    \filldraw[fill=black!50, draw=black,opacity=0.1] (0,1.5) -- (0,3.5) -- (1.5,4.5) -- (3,4.5) -- (4.5,3) -- (4.5,1.5) -- (3,0) -- (1.5,0) -- cycle; 
    \end{scope}
    \end{scope}

    \draw[style=help lines, line width = 0.01mm,lightgray!60] (-0.5,-0.5) grid[step=0.5cm] (5,5); 
    
    \filldraw[fill=blue!30, draw=black, opacity = 0.6] (0,3.5) --  (1.5,4.5) -- (2.5,4.5) -- cycle;

    \foreach \x in {0,0.5,...,4.5}{                           
        \foreach \y in {0,0.5,...,4.5}{                       
        \node[draw,circle,inner sep=0.01pt,fill,lightgray!60] at (\x,\y) {}; 
        }
    }
    \node[above left] (x) at (1.5,4.5) {$\hat\bx$};
    \node[draw,circle,inner sep=1.2pt,fill,darkgray] at (1.5,4.5) {};
    
    \node[above right] (x) at (3,4.5) {$\bx^7$};
    \node[draw,circle,inner sep=0.8pt,fill,darkgray] at (3,4.5) {};
    
    \node[right] (x) at (4.5,3) {$\bx^6$};
    \node[draw,circle,inner sep=0.8pt,fill,darkgray] at (4.5,3) {};
    
    \node[right] (x) at (4.5,1.5) {$\bx^5$};
    \node[draw,circle,inner sep=0.8pt,fill,darkgray] at (4.5,1.5) {};
    
    \node[below] (x) at (3,0) {$\bx^4$};
    \node[draw,circle,inner sep=0.8pt,fill,darkgray] at (3,0) {};

    \node[below] (x) at (1.5,0) {$\bx^3$};
    \node[draw,circle,inner sep=0.8pt,fill,darkgray] at (1.5,0) {};

    \node[below] (x) at (0,1.5) {$\bx^2$};
    \node[draw,circle,inner sep=0.8pt,fill,darkgray] at (0,1.5) {};

    \node[above] (x) at (0,3.5) {$\bx^1$};
    \node[draw,circle,inner sep=0.8pt,fill,darkgray] at (0,3.5) {};
    
    \node[below left] (x) at (-0.5,-0.5) {$(0,0)$};
    \node[draw,circle,inner sep=0.8pt,fill,darkgray] at (-0.5,-0.5) {};
    
    \node[below] (x) at (4.5,-0.5) {$(10,0)$};
    \node[draw,circle,inner sep=0.8pt,fill,darkgray] at (4.5,-0.5) {};

    \node[left] (x) at (-0.5,4.5) {$(0,10)$};
    \node[draw,circle,inner sep=0.8pt,fill,darkgray] at (-0.5,4.5) {};
    \node[draw,diamond,inner sep=1.3pt,fill,blue!60] at (1,3.5) {};
    
    \node[draw,diamond,inner sep = 1.3pt,fill,blue!60] at (1.5,3) {};
    
    \node[draw,diamond,inner sep = 1.3pt,fill,blue!60] at (2,3) {};

    \node[draw,diamond,inner sep = 1.3pt,fill,blue!60] at (2.5,3.5) {};

    \node[draw,diamond,inner sep = 1.3pt,fill,blue!60] at (2.5,4) {};

    \node[draw,diamond,inner sep = 1.3pt,fill,blue!60] at (2.5,4.5) {};

    
    \fill[fill=red!55] (1.5,4.5) -- (1.8,4.5) arc (0:-145:0.3cm) -- cycle;
    
    \draw[red!50, line width = 0.4mm, dashed] (1.5,4.5) circle (0.3cm);

    \draw [->, line width = 0.05mm,dashed,darkgray] (1.5,4.5) -- (0,3.5);
    \draw [->, line width = 0.3mm,darkgray] (1.5,4.5) -- (0.05,3.51);
    
    \draw [->, line width = 0.05mm,dashed,darkgray] (1.5,4.5) -- (0,1.5);
    \draw [->, line width = 0.3mm, darkgray] (1.5,4.5) -- (1.1,3.7);
    
    \draw [->, line width = 0.05mm,dashed,darkgray] (1.5,4.5) -- (1.5,0);
    \draw [->, line width = 0.3mm,darkgray] (1.5,4.5) -- (1.5,3.2);

    \draw [->, line width = 0.05mm,dashed,darkgray] (1.5,4.5) -- (3,0);
    \draw [->, line width = 0.3mm,darkgray] (1.5,4.5) -- (1.9,3.3);
    
    \draw [->, line width = 0.05mm,dashed,darkgray] (1.5,4.5) -- (4.5,1.5);
    \draw [->, line width = 0.3mm,darkgray] (1.5,4.5) -- (2.25,3.75);
    
    \draw [->, line width = 0.05mm,dashed,darkgray] (1.5,4.5) -- (4.5,3);
    \draw [->, line width = 0.3mm,darkgray] (1.5,4.5) -- (2.25,4.125);
    
    \draw [->, line width = 0.05mm,dashed,darkgray] (1.5,4.5) -- (3,4.5);
    \draw [->, line width = 0.3mm,darkgray] (1.5,4.5) -- (2.3,4.5);

\draw[->, thick] (-0.5,-0.5)--(5,-0.5) node[below]{};
\draw[->, thick] (-0.5,-0.5)--(-0.5,5) node[above]{};

    \end{scope}
\end{tikzpicture}}
\vspace{-0.3cm}
\caption{$\mX \subset \mathbb{Z}^n$}\label{fig:exp_con_ball}
\end{subfigure}
\hspace{-0.5cm}
\begin{subfigure}[t]{0.48\textwidth}
\centering
\scalebox{0.97}{
\begin{tikzpicture}
    \centering
    
    \filldraw[fill=blue!70, draw=blue!70, opacity = 0.15] (-1.55,2.38) --  (-0.01,2.38) -- (-0.01,0) -- cycle; 
    
    \draw[line width=0mm,   black!40] (-2.4,0) -- (2.4,0) node [below] {};
    \draw[line width=0.1mm,   blue!70] (-0.04,2.4) -- (-0.04,-2.4) node [below left] {\color{black!70}$\bx^7$};
    \draw[->,line width = 0.3mm, blue!60] (-0.04, 2.1) -- (-0.35,2.1);
    \draw[->,line width = 0.3mm, blue!60] (-0.04, -2.1) -- (-0.35,-2.1);
    \draw[line width=0.1mm,   blue!70] (-1.6,2.4) node [above] {\color{black!70}$\bx^1$} -- (1.6,-2.4);
    \draw[->,line width = 0.3mm, blue!60] (-1.2, 1.8) -- (-0.92,1.98);
    \draw[->,line width = 0.3mm, blue!60] (1.2, -1.8) -- (1.48,-1.62);
    \draw[line width=0mm,   blue!70] (-2.4,1.2) node [left] {\color{black!70}$\bx^2$} -- (2.4,-1.2);
    \draw[->,line width = 0.3mm, blue!60] (-1.8, 0.9) -- (-1.65,1.2);
    \draw[->,line width = 0.3mm, blue!60] (1.8, -0.9) -- (1.95,-0.6);
    \draw[line width=0mm,   blue!70] (-2.4,-0.8) node [left] {\color{black!70}$\bx^4$} -- (2.4,0.8);
    \draw[->,line width = 0.3mm, blue!60] (-1.93, -0.65) -- (-2.02,-0.35);
    \draw[->,line width = 0.3mm, blue!60] (1.93, 0.65) -- (1.84,0.95);
    \draw[line width=0mm,   blue!70] (-2.4,-2.4) node [below] {\color{black!70}$\bx^5$} -- (2.4,2.4);
    \draw[->,line width = 0.3mm, blue!60] (-1.6, -1.6) -- (-1.82,-1.38);
    \draw[->,line width = 0.3mm, blue!60] (1.6, 1.6) -- (1.38,1.82);
    \draw[line width=0mm,   blue!70] (-1.2,-2.4) node [below] {\color{black!70}$\bx^6$} -- (1.2,2.4);
    \draw[->,line width = 0.3mm, blue!60] (-1,-2) -- (-1.25,-1.875);
    \draw[->,line width = 0.3mm, blue!60] (1, 2) -- (0.75,2.125);
    \draw[line width=0mm,   blue!70] (-2.4,0) node [above left] {\color{black!70}$\bx^3$} -- (2.4,0);
    \draw[->,line width = 0.3mm, blue!60] (-2,0) -- (-2,0.3);
    \draw[->,line width = 0.3mm, blue!60] (2,0) -- (2,0.3);

    \draw[thick,<->] (-2.6,0)--(2.6,0) node[right] {$\bc_1$}; 
    \draw[thick,<->] (0,-2.6)--(0,2.6) node[above] {$\bc_2$}; 

\end{tikzpicture}}
\vspace{-0.3cm}
\caption{$\mC(\bhx,\mE(\mX))$}\label{fig:half_spaces}
\end{subfigure}
\vspace{0.1cm}
\caption{(a) The convex hull of a forward-feasible region $\mX \subset \mathbb{Z}^2$ defined on a two-dimensional lattice is shaded with the points in $\mE(\mX)$ marked. (b) The half-space in the cost vector space defined by each extreme point, and the inverse-feasible region lying at the intersection of half-spaces shaded.}\label{fig:main_example_fig} 
\end{figure}

These observations are illustrated in Figure~\ref{fig:main_example_fig}. The extreme points $\{\bx^1,...,\bx^7\}$ in Figure~\ref{fig:exp_con_ball} induce half spaces for $\bc \in \mathbb{R}^n$ and the intersection of these half-spaces is the cone shown in Figure~\ref{fig:half_spaces}. However, the same cone can be formed using the (blue) diamond points, all of which are \emph{interior} forward-feasible points, plus $\bx^1$, in Figure~\ref{fig:exp_con_ball}. 

The last observation above is formalized in the following lemma. Let $\by(\lambda, \bhx, \bx) :=  \hat\bx + \lambda(\bx - \hat\bx)$, $\lambda > 0$, define a point along the ray from $\hat\bx$ to $\bx$.

\begin{lemma}\label{lemma:rays_equiv}
Let $\mE(\mX) = \{\bx^1,\ldots,\bx^N\}$ and $\bar{\mE}(\mX) = \{\by(\lambda_1,\bhx,\bx^1), \ldots , \by(\lambda_N, \bhx,\bx^N)\}$ for arbitrary positive scalars $\lambda_i > 0, i = 1, \ldots, N$. Then,
\[\mC(\hat\bx, \mE(\mX)) = \mC(\hat\bx, \bar{\mE}(\mX)). \]
\end{lemma}


Lemma \ref{lemma:rays_equiv} shows that the set of extreme points $\mE(\mX)$ is not necessary to characterize $\mC(\hat\bx,\mX)$ and can instead be replaced with a set of interior and/or exterior points. This observation can be generalized with a necessary and sufficient condition for any collection of points to be a generator set. Let $\mY(\bhx, \mX)$ denote the polyhedral cone pointed at $\hat\bx$ generated by $\{\by(1,\bhx,\bx)\}_{\bx\in \mX}$.

\begin{theorem}\label{thm:certificate_set} A set $\mG \subseteq \mathbb{R}^n$ is a generator set if and only if $\mY(\bhx,\mX) = \mY(\bhx,\mG)$. 
\end{theorem}

While Lemma \ref{lemma:rays_equiv} shows that generator sets can be constructed without using any extreme points, Theorem \ref{thm:certificate_set} goes further and shows that far fewer points than the number of extreme points may be sufficient to construct a generator set. As long as the set of extreme rays of $\mY(\bhx, \mG)$ are the same as those of $\mY(\bhx, \mX)$, then $\mG$ is a generator set -- all other points in $\mG$ that do not form an extreme ray are redundant in defining the inverse-feasible region. Referring back to Figure \ref{fig:exp_con_ball}, this result implies that any two points $\by(\lambda_1,\bhx,\bx^1)$ and $\by(\lambda_7, \bhx,\bx^7)$ with $ \lambda_1, \lambda_7 > 0$ are sufficient to define the generator set in that example. In Figure \ref{fig:half_spaces}, the inverse-feasible region is exactly the intersection of the half-spaces determined by $\bhx-\bx^1$ and $\bhx - \bx^7$.

Note that for any forward-feasible point $\bx \in \mX$, the ray $\{\by(\lambda,\bhx,\bx)\}_{ \lambda \geq 0}$ is in $\mY(\bhx,\mX)$. Thus, $\mY(\bhx, \mG) \subseteq \mY(\bhx, \mX)$ for any collection of forward-feasible points $\mG \subseteq \mX$. Then, the set $\mG \subseteq \mX$ is a forward-feasible generator set if $\mY(\bhx,\mX) \subseteq \mY(\bhx,\mG)$, or equivalently if $\mX \subseteq \mY(\bhx,\mG)$. This is shown by the following corollary.

\begin{corollary}\label{cor:ff_cert_set}
A set $\mG \subseteq \mathbb{R}^n$ is a forward-feasible generator set if and only if $\mG \subseteq \mX \subseteq \mY(\bhx, \mG)$.
\end{corollary} 

Analyzing forward-feasible generator sets can be particularly informative for the design of cutting plane algorithms. For example, inverse MILO problems can be solved by the cutting plane algorithm presented in \cite{wang2009cutting} because $\mE(\mX)$ is a forward-feasible generator set. Below, we present 
an equivalent characterization of forward-feasible generator sets that does not rely on the notion of rays. This result provides useful intuition to aid in designing a new cutting plane algorithm. Let $\mB(\epsilon,\hat\bx)$ denote a closed ball of radius $\epsilon > 0$ around $\hat\bx$. 

\begin{theorem}\label{thm:sufficient_certificates} A set $\mG \subseteq \mathbb{R}^n$ is a forward-feasible generator set if and only if $\mG \subseteq \mX$ and there exists an $\epsilon>0$ such that 
\begin{align}
\mB(\epsilon,\hat\bx) \cap \conv(\mX) \subseteq \conv(\mG \cup \{\hat\bx\}).
\end{align}
\end{theorem}
Theorem~\ref{thm:sufficient_certificates} states that a set $\mG \subseteq \mX$ is a forward-feasible generator set if and only if the convex hull of the set of points in $\mG$ and $\{\hat\bx\}$ \emph{contains} an epsilon ball around $\hat\bx$ intersected with the convex hull of $\mX$. Going back to Figure \ref{fig:exp_con_ball}, $\mG = \{\bx^1, \bx^7\}$ is a forward-feasible generator set since the convex hull of $\{\bx^1, \bx^7, \bhx\}$ contains the part of an epsilon ball around $\bhx$ inside $\conv(\mX)$, which is shaded in red. Considering the ray $\{\by(\lambda_1, \bhx,\bx^1)\}_{\lambda_1 \geq 0}$, only $\bx^1$ is a forward-feasible point, which means that any forward-feasible generator set must include $\bx^1$. However, there are multiple forward-feasible points along the ray $\{\by(\lambda_7, \bhx,\bx^7)\}_{\lambda_7 \geq 0}$, which means that any of those feasible points, namely (5,10), (6,10) or (7,10), can be part of the forward-feasible generator set.

The two key takeaways from Theorem~\ref{thm:sufficient_certificates} are that a forward-feasible generator set can: (i) consist of many fewer than $|\mE(\mX)|$ points, and (ii) be generated from non-extreme points that lie much closer to $\bhx$.


\section{Cutting Plane Algorithm with Trust Regions}\label{sec:general_framework}

The theoretical insights from the previous section indicate that a new approach to generating cuts may yield a more effective solution algorithm. 
Instead of generating cuts by identifying extreme points of the convex hull of the forward-feasible region $\mX$, we generate cuts by identifying interior points of $\mX$, which we accomplish using \emph{trust regions}. More specifically, we intersect $\mX$ with a trust region, and generate cuts at interior points of $\mX$ obtained as extreme points of the convex hull of this restriction of $\mX$.  

Let $\mT(\hat\bx,p)$ be a trust region around the point $\hat\bx$ of size $p \geq 1$, defined as
\begin{equation}\label{eq:normed_TR}
\mT(\bhx, p) := \{\by \in \mathbb{R}^n \ | \ \norm{\bhx - \by}_1 \leq p \}.
\end{equation}
We simplify notation to $\mT$ when the parameter $p$ is not the focus of discussion.

Guided by Theorem \ref{thm:sufficient_certificates}, we define trust regions to be \emph{centered at $\bhx$}. The trust region, when imposed onto $\mX$, represents an attempt to capture a forward-feasible generator set. The 1-norm in equation \eqref{eq:normed_TR} is chosen because the trust region has a polynomial number of extreme points. When intersected with $\mX$, the resulting restricted forward-feasible region $\mT \cap \mX$ itself is more likely to have a smaller number of extreme points compared to $\mX$. This can reduce the number of forward-feasible points that need to be identified in order to build a generator set. Furthermore, the extreme points of this smaller, restricted region may also be computationally easier to identify. These considerations will be discussed in detail below.

\subsection{Model decomposition and cutting plane framework}\label{subsec:trustregions}

The inverse MILO problem is decomposed into a master problem, describing a relaxation of the inverse model \eqref{model:IOM-Bilevel1}, and a cut generation subroutine that is used to iteratively tighten the master problem. The framework of our cutting plane algorithm is given in Algorithm \ref{algo:CP-TR}. A global information set $\mathscr{I}$ is used to pass information between successive iterations of the subroutine, which specifically includes the trust region and outer-loop index for our purposes. It is initialized with a trust region of size $p^0$ centered at $\hat\bx$, and the iteration index $i=0$.

\begin{algorithm}[h]
\caption{A general cutting plane algorithm for inverse MILO}
\label{algo:CP-TR}
\vspace{0.3cm}
\textbf{Input:} An inverse MILO problem instance $(\bc^0, \hat\bx, \mX)$, initial trust region size $p^0$\\
\textbf{Output:} An inverse-optimal solution $\bc^*$
\begin{algorithmic}[1]
\State Initialize $i=0, \, \mathscr{I}^0 = (\mT(\hat\bx,p^0), i), \, \tilde\mX^i = \emptyset, \, \tilde\bc^i = \bc^0$ 
\State Run \textbf{SUBROUTINE}($\tilde\bc^i, \bhx, \mX, \mathscr{I}^{i})$, let $\tilde\bx^{i}$ and $\mathscr{I}^{i+1}$ be its output
\While{$\tilde\bx^i \neq \bhx$}
\State $i \gets i+1$
\State $\tilde\mX^{i} \gets \tilde\mX^{i-1} \cup \{\tilde\bx^{i-1}\}$
\State Solve $\textbf{MP}(\bhx,\tilde\mX^{i})$, let $\tilde\bc^i$ be its optimal solution \label{step:candidate_cost}
\State Run \textbf{SUBROUTINE}($\tilde\bc^i, \bhx, \mX, \mathscr{I}^{i})$, let $\tilde\bx^{i}$ and $\mathscr{I}^{i+1}$ be its output
\EndWhile \\
\Return $\bc^* = \tilde\bc^i$ 
\end{algorithmic}
\end{algorithm}

The master problem is defined below, where $\tilde\mX \subseteq \mX$ denotes a finite set of forward-feasible points. 
\begin{subequations}
\label{model:MP}
\begin{align}
\textbf{MP}(\hat\bx, \tilde\mX): \ \ \underset{\iz{\bc \in \mP}}{\text{minimize}} \quad & \norm{\bc - \bc^0}_1 \\
\text{subject to} \quad & \bc^\top (\hat\bx - \bx) \leq 0, \quad \forall \bx \in \tilde\mX,\label{cons:decomp_master}\\
& \by^\top \bA = \bc^\top,\label{cons:duality1}\\
& \by \geq \textbf{0}.\label{cons:duality2}
\end{align}
\end{subequations}
Constraints \eqref{cons:decomp_master} ensure that $\bhx$ is optimal with respect to $\bc$ over $\tilde\mX$. As mentioned in Remark \ref{rmk:Cbdd}, the inverse-feasible region $\mC(\bhx,\mX)$ only contains cost vectors for which the forward optimization problem is bounded. However, due to its relaxed nature, the master problem can propose cost vectors that lead to unbounded forward problems in the subroutine. This can be prevented by constraints \eqref{cons:duality1} and \eqref{cons:duality2}, as proposed in \cite{wang2009cutting}, which are not needed if the forward-feasible region is bounded. Let $\tilde\bc$ denote an optimal solution to the master problem $\textbf{MP}(\hat\bx, \tilde\mX)$. 

\subsection{Cut generation subroutine}\label{subsec:CGS}

Given $\tilde\bc$, the cut generation subroutine either finds a feasible point $\tilde\bx \in \mX$ satisfying $\tilde\bc^\top(\hat\bx - \tilde\bx) > 0$ (i.e.,  $\tilde\bc \notin \mC(\hat\bx, \tilde\mX \cup \{\tilde\bx\})$), which generates a cut, or returns $\bhx$, verifying that $\tilde\bc \in \mC(\bhx,\mX)$. Note that the verification in the latter case requires solving the forward problem $\FP(\tilde\bc,\mX)$, which means the cut generation subroutine is at least as hard.  

The cut generation subroutine is presented in Subroutine \ref{subroutine:CP-TR}. It consists of four key components: the cut generation problem (Step \ref{step:solve_ropt}) and three functions to modify the trust region (Steps \ref{step:REMOVEfunction}, \ref{step:SAVEfunction} and \ref{step:updateTR}), which we describe next.

\begin{subroutine}[ht]
\caption{Cut generation subroutine}
\label{subroutine:CP-TR}
\vspace{0.3cm}
\textbf{Input:} A candidate objective $\tilde \bc$, forward-feasible point $\bhx$, forward-feasible region $\mX$, information set $\mathscr{I}^{\textrm{in}} = \{\mT, i\}$\\ 
\textbf{Output:} A forward-feasible point $\tilde \bx$ and updated information set $\mathscr{I}^{\textrm{out}}$
\begin{algorithmic}[1]
\State Initialization: $k = 1$, $\mT^0 = \mT$. \label{step:initialization}
\State $\mT^k \gets \text{REMOVE}(\mT^{k-1})$ \label{step:REMOVEfunction}
\State Solve $\textbf{FP}(\tilde\bc, \mT^k \cap \mX)$, let $\tilde\bx$ be the solution \label{step:solve_ropt} 
\If {$\tilde\bc^\top(\bhx-\tilde\bx) > 0$} 
    \State $\mT^* \gets \text{SAVE}(\mT^k)$,
    \label{step:SAVEfunction}
    \State \textbf{Return} $\tilde\bx$, $\mathscr{I}^{\textrm{out}} = \{\mT^*, i+1\}$\label{step:return_after_SAVE}
\ElsIf{$\mT^k \neq \mathbb{R}^n$} \label{step:SP_end_loop}
    \State $\mT^{k+1} \gets $ UPDATE($\mT^k$)     \label{step:updateTR} 
    \State $k \gets k+1$\label{step:update_k}
    \State Return to Step \ref{step:REMOVEfunction} \label{step:loop}
\Else \label{step:terminate}
    \State \textbf{Return} $\bhx$, $\mathscr{I}^{\textrm{in}}$ \label{step:finalstep}
\EndIf
\end{algorithmic}
\end{subroutine}

\subsubsection{Cut generation problem (Step \ref{step:solve_ropt}).} 
Given a candidate objective $\tilde\bc$ and a trust region $\mT$, the cut generation problem is the forward problem solved over $\mT \cap \mX$ \iz{(i.e., $\FP(\tilde\bc,\mT \cap \mX)$)}, which we refer to as a \emph{subregion} of $\mX$. We note that the choice of 1-norm in the definition of employed trust regions ensures that the cut generation problem remains linear. 

If the optimal solution of the cut generation problem $\tilde\bx$ satisfies the  condition $\tilde\bc^\top(\hat\bx - \tilde\bx) > 0$, then $\tilde\bx$ is returned; we call such an $\tilde\bx$ a \emph{violated forward-feasible point}. If no cut is found, then there are two possibilities depending on whether a trust region was used. If no trust region was used, i.e., $\mT^k = \mathbb{R}^n$, then $\bhx$ is returned and Algorithm \ref{algo:CP-TR} terminates with $\tilde\bc$. Otherwise, $\tilde\bc$ renders $\bhx$ optimal over a subregion of $\mX$, but not necessarily over $\mX$. Thus, larger subregions need to be considered before it can be verified that no additional violated cuts exist, which leads to a trust region update. 

\subsubsection{UPDATE function (Step \ref{step:updateTR}).} 
When a subregion $\mT(\bhx, p) \cap \mX$ has been exhausted of violated cuts, we increase the size $p$ of the trust region to $\delta p$ for some $\delta > 1$:
\begin{equation}\label{eq:update-centered}
\text{UPDATE}(\mT(\bhx, p)) = \mT(\bhx, \delta p). 
\end{equation}
For example, setting $\delta = 2$ doubles the size of the trust region each time UPDATE is called. There is a natural trade-off to be considered when choosing $\delta$. If $\delta$ is too small, then the trust regions grow slowly and may lead to the generation of subregions with no violated cuts. This will result in time-consuming, ``empty" calls to the cut generation problem, i.e., calls in which the cut generation problem is solved but is unable to generate any violated forward-feasible point. On the other hand, if $\delta$ is too large, the computational advantages of solving the cut generation problem over a smaller subregion with fewer extreme points may quickly be lost.

\subsubsection{REMOVE function (Step \ref{step:REMOVEfunction}).} The REMOVE function periodically removes the existing trust region (by setting $\mT = \mathbb{R}^n$) at regular intervals in both outer (i.e., master) iterations $i$ and inner (i.e., subroutine) iterations $k$. In particular, for fixed index values $i^* > 0$ and $k^* > 0$,
\begin{equation}\label{eq:eval_function} 
\text{REMOVE}(\mT) = \begin{cases}
\mathbb{R}^n, & \text{ if } i \in \{i^*, 2i^*, \ldots\} \text{ or } k = k^*,\\ 
\mT, & \text{ otherwise. }
\end{cases}
\end{equation}
In other words, the cut generation problem is solved over $\mX$ every $i^*$ iterations in the outer loop and when the inner loop reaches the iteration limit of $k^*$. The REMOVE function guarantees finite convergence of Algorithm~\ref{algo:CP-TR}, irrespective of how the UPDATE function is defined, because an extreme point of $\mX$ is returned at least every $i^*$ iterations in a finite number of (at most $k^*$) attempts to verify that a candidate objective $\tilde\bc^i$ is inverse-feasible when no more violated cuts exist.

Note that the classical cutting plane algorithm proposed in \cite{wang2009cutting} is a special case of our trust region-based approach where $i^* = k^* = 1$. 

\subsubsection{SAVE function (Step \ref{step:SAVEfunction}).} 

The SAVE function is called when a forward-feasible point yielding a violated cut is found by the cut generation problem. If such a point is found during an iteration where the trust region is $\mathbb{R}^n$ (i.e., when REMOVE is called), the SAVE function saves the previous trust region. On the other hand, if the REMOVE function did not remove the trust region during the current iteration, then the current trust region is saved. This ensures that there always exists a ``non-trivial" trust region ($\mT(\bhx,p) < \infty$) saved.

Mathematically, the SAVE function is defined as:
\begin{equation*} 
\text{SAVE}(\mT^k) = \begin{cases}
\mT^{k-1} & \text{ if } \mT^{k}  = \mathbb{R}^n,\\ 
\mT^k & \text{ otherwise. } 
\end{cases}
\end{equation*}

\subsubsection{Discussion and example.}
Trust regions can reduce the difficulty of solving the cut generation problem by reducing the size of the forward-feasible region. However, they can also lead to the generation of \emph{stronger} cuts.%

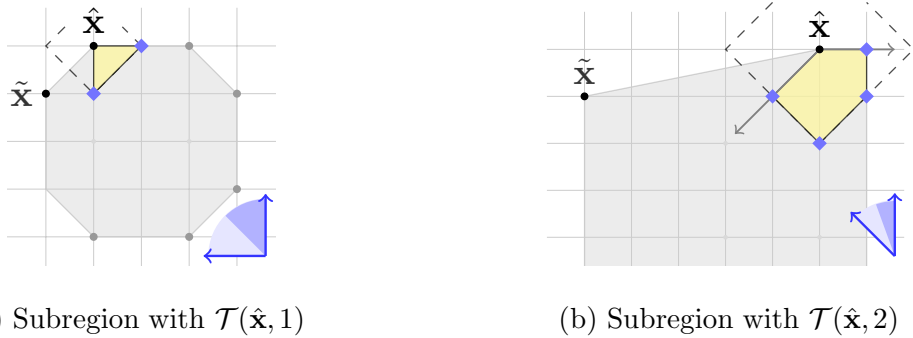
\begin{figure}[h]
\centering
\begin{subfigure}[t]{0.47\textwidth}
\centering
\scalebox{1.27}{
  \begin{tikzpicture}
    \begin{scope}
    \clip (0.6,2.2) rectangle (3.4cm,5.1cm); 
    \begin{scope}[transparency group]
    \begin{scope}[blend mode=multiply]
    \filldraw[fill=black!50, draw=black,opacity=0.15] (1,3) -- (1,4) -- (1.5,4.5) -- (2.5,4.5) -- (3,4) -- (3,3) -- (2.5,2.5) -- (1.5,2.5) -- cycle; 

    \end{scope}
    \end{scope}

    \draw[style=help lines, line width = 0.1mm,lightgray!80] (0.5,1.48) grid[step=0.5cm] (3.5,4.9); 
    
    \draw[draw=black, dashed, opacity = 0.7] (1.2,4.7) -- (1,4.5) --  (1.5,4) -- (2,4.5) -- (1.8,4.7);
    \filldraw[fill=yellow!50, opacity = 0.7] (1.5,4) -- (1.5,4.5) -- (2,4.5) -- cycle;
     \fill[fill=blue!10] (3.3,2.3) -- (2.7,2.3) arc (180:90:0.6cm) -- cycle;  
     \fill[fill=blue!30] (3.3,2.3) -- (2.88,2.72) arc (135:90:0.6cm) -- cycle;  
    \draw [->, line width = 0.25mm,blue!80] (3.3,2.3) -- (2.65,2.3);
    \draw [->, line width = 0.25mm,blue!80] (3.3,2.3) -- (3.3,2.95);
    
    \foreach \x in {1.5,...,3.5}{                           
        \foreach \y in {1.5,...,3.5}{                       
        \node[draw,circle,inner sep=0.01pt,fill,lightgray!60] at (\x,\y) {}; 
        }
    }
    \node[draw,diamond,inner sep=1pt,fill,blue!55] at (1.5,4) {};
    \node[draw,diamond,inner sep=1pt,fill,blue!55] at (2,4.5) {};
    \node[draw,circle,inner sep=0.7pt,fill,black] at (1.5,4.5) {};
    \node[draw,circle,inner sep=0.7pt,fill,gray!80] at (2.5,4.5) {};
    \node[draw,circle,inner sep=0.7pt,fill,gray!80] at (3,4) {};
    \node[draw,circle,inner sep=0.7pt,fill,gray!80] at (3,3) {};
    \node[draw,circle,inner sep=0.7pt,fill,gray!80] at (2.5,2.5) {};
    \node[draw,circle,inner sep=0.7pt,fill,gray!80] at (1.5,2.5) {};
    \node[above] (xhat) at (1.5,4.5) {$\bhx$};
    \node[left] (xn) at (1,4) {\color{black!80}$\tilde\bx$};
    \node[draw,circle,inner sep=0.7pt,fill,black] at (1,4) {};
    \end{scope}
\end{tikzpicture}}
\caption{Subregion with $\mT(\bhx, 1)$}\label{fig:one_manhattan_norm}
\end{subfigure}
\begin{subfigure}[t]{0.47\textwidth}
\centering
\scalebox{1.25}{
  \begin{tikzpicture}
    \begin{scope}
    \clip (-0.4,2.2) rectangle (3.6cm,5.2cm); 
    \begin{scope}[transparency group]
    \begin{scope}[blend mode=multiply]
    \filldraw[fill=black!50, draw=black,opacity=0.15] (0,2) -- (0,4) -- (2.5,4.5) -- (3,4.5) -- (3,3.5) -- (3,2) -- (2.5,1.5) -- (1.5,1.5) -- cycle; 

    \end{scope}
    \end{scope}

    \draw[style=help lines, line width = 0.1mm,lightgray!80] (-0.4,1.95) grid[step=0.5cm] (3.4,4.9); 
    
    \draw[draw=black, dashed, opacity = 0.7] (2,5) -- (1.5,4.5) -- (2.5, 3.5) -- (3.5,4.5) -- (3,5);
    \filldraw[fill=yellow!50, opacity = 0.7] (2,4) -- (2.5,4.5) -- (3,4.5) -- (3,4) -- (2.5, 3.5) -- cycle;

    \draw [->, line width = 0.2mm,gray] (2.5,4.5) -- (1.6,3.6);
    \draw [->, line width = 0.2mm,gray] (2.5,4.5) -- (3.3,4.5);
    \foreach \x in {1.5,...,3}{                           
        \foreach \y in {1.5,...,3.5}{                       
        \node[draw,circle,inner sep=0.01pt,fill,lightgray!60] at (\x,\y) {}; 
        }
    }
    \node[draw,diamond,inner sep=1pt,fill,blue!55] at (2,4) {};
    \node[draw,diamond,inner sep=1pt,fill,blue!55] at (3,4) {};
    \node[draw,diamond,inner sep=1pt,fill,blue!55] at (3,4.5) {};
    \node[draw,diamond,inner sep=1pt,fill,blue!55] at (2.5,3.5) {};
    \node[draw,circle,inner sep=0.7pt,fill,black] at (2.5,4.5) {}; 
    \node[above] (x2) at (0,4) {\color{black!80}$\tilde\bx$};
    \node[above] (bhx) at (2.5,4.5) {$\bhx$};
    \node[draw,circle,inner sep=0.7pt,fill,black] at (0,4) {};
     \fill[fill=blue!10] (3.3,2.3) -- (2.88,2.72) arc (135:90:0.6cm) -- cycle;  
     \fill[fill=blue!30] (3.3,2.3) -- (3.1,2.85) arc (110:90:0.6cm) -- cycle;  
    \draw [->, line width = 0.25mm,blue!80] (3.3,2.3) -- (2.83,2.77);
    \draw [->, line width = 0.25mm,blue!80] (3.3,2.3) -- (3.3,2.97);
    
    \end{scope}
\end{tikzpicture}}
\caption{Subregion with $\mT(\bhx,2)$}\label{fig:one_infinity_norm_extreme}
\end{subfigure}
\vspace{0.3cm}
\caption{Two examples of when a trust region is applied on a two-dimensional integer feasible region. The convex hull of $\mX$ and $\mX \cap \mT$ is shaded in gray and yellow, respectively. The shaded cone at the right bottom corner of each figure represents $-\mC(\bhx,\mX\cap\mT)$ while the darker shaded cone represents $-\mC(\bhx,\mX)$, a subset of $-\mC(\bhx,\mX\cap\mT)$.}\label{fig:approximate_certificate_set}
\end{figure}

We use Figure \ref{fig:one_manhattan_norm} to illustrate this latter point. A trust region of size one ($p = 1$) is imposed on a two-dimensional integer forward-feasible region, creating a subregion containing two extreme points (indicated using blue diamonds), excluding $\bhx$. Those two extreme points of the subregion have the same cut generation potential as the five extreme points of $\mX$ labeled in gray. More specifically, the set of cuts generated by the former set of extreme points can dominate those generated by the latter. Once this trust region is exhausted, it must be updated. 
In Figure \ref{fig:one_manhattan_norm}, increasing the size of the trust region by a single unit would cover $\tilde\bx$, the only additional point needed to build a generator set. The trust region approach guarantees that at most three cuts need to be generated, or equivalently three forward-feasible points need to be identified, to solve the inverse MILO problem. On the other hand, the cutting plane algorithm without trust regions, which computes extreme points of $\mX$, may require generating many more cuts.

The REMOVE function periodically removes trust regions and provides opportunities for the algorithm to compute points outside of the existing trust region. We consider another example in Figure \ref{fig:one_infinity_norm_extreme} with a trust region of $\mT(\bhx,2)$ over a different two-dimensional integer forward-feasible region. Assuming the subregion $\mX \cap \mT(\bhx,2)$ has been exhausted of cuts, removal of this trust region would lead to the computation of $\tilde\bx$ with certainty, as it is the only extreme point of $\conv(\mX)$ that can still generate a valid cut. This would terminate the algorithm at the next step. On the other hand, if the size of the trust region was increased incrementally each time without removing the trust region, many unnecessary interior points would be computed before a sufficiently large trust region captures $\tilde\bx$. Periodic removal of trust regions thus helps to reduce the potential of being overly-conservative, and when used in conjunction with trust regions create an effective cut generation subroutine.

\section{Algorithm Enhancements} \label{sec:enhancements}

In this section, we propose two enhancements to Subroutine \ref{subroutine:CP-TR} that can be implemented independently or together. These two enhancements are introduced to improve cut generation speed, particularly during iterations when the subregion is large. 

\subsection{Early-stop heuristic for the cut generation problem}\label{subsec:early-stop}
As defined in Subroutine \ref{subroutine:CP-TR}, the cut generation problem is solved to optimality in each iteration. However, this can be time consuming, particularly during iterations where the trust region is large or removed entirely. Here, we introduce a solution time threshold $\tau$. If the cut generation problem is not solved by $\tau$, but at least one violated forward-feasible point has been found, the feasible point of maximum violation is returned. If no violated forward-feasible point has been found by $\tau$, the cut generation problem will return the first such point found after $\tau$. This enhancement can reduce the cut generation time and can be effective for two reasons: (i) Solving the cut generation problem may not necessarily provide the best cuts anyways, as discussed in Section \ref{sec:backgroundandmotivation}, and (ii) the solver may have already found an optimal violated forward-feasible point in the cut generation problem but has not computed a bound strong enough to prove its optimality. When a standard MILO solver, such as CPLEX and Gurobi, is used to solve the forward problem in Subroutine \ref{subroutine:CP-TR}, the early-stop heuristic can be implemented using a callback function. 

\subsection{Stochastic dimensionality reduction of the trust region}

This enhancement is motivated by the observation that a violated forward-feasible point generated within a trust region may have many components with the same value as the corresponding components in $\bhx$. For example, if $\mX$ contains only integer variables, then all forward-feasible points within a trust region $\mathcal{T}(\hat\bx,p)$ around $\hat\bx$ will differ from $\hat\bx$ by at most $p$ components. By explicitly restricting \emph{which} subset of components can differ from the values of $\bhx$, we can decrease the dimensionality and hence the size of the cut generation problem.

Let $\mS \subset \mN := \{1,\hdots,n\}$ and 
\begin{equation}
    \mT_{\mS}(\bhx,p) := \{\by \in \mathbb{R}^n \ | \ \by \in \mT(\bhx, p), \ y_j = \hat{x}_j \ \forall j \notin \mS\} \label{eq:TR_with_DR}
\end{equation}
define a lower dimensional trust region containing points that can differ from $\hat\bx$ only in the indices contained in $\mS$. 

\begin{subroutine}[ht]
\caption{Cut generation subroutine with dimensionality reduction}
\label{subroutine:CP-TRdimred}
\vspace{0.3cm}
\textbf{Input:} A candidate objective $\tilde \bc$, forward-feasible point $\bhx$, collection of points $\tilde\mX$, information set $\mathscr{I}^{\textrm{in}} = \{\mT, i, h\}$\\ 
\textbf{Output:} A forward-feasible point $\tilde \bx$ and updated information set $\mathscr{I}^{\textrm{out}}$
\begin{algorithmic}[1]
\State Initialization: $k = 1$, $\mT^0 = \mT$. 
\State $\mT^k \gets \text{REMOVE}(\mT^{k-1})$ 
\State Solve $\textbf{FP}(\tilde\bc, \mT^k \cap \mX)$, let $\tilde\bx$ be the solution 
\If {$\tilde\bc^\top(\bhx-\tilde\bx) > 0$} 
    \State $\mT^* \gets \text{SAVE}(\mT^k)$
     \If {$\mT^k \neq \mathbb{R}^n$}
     \State $h \gets 0$ \label{step2:h_reset}
     \EndIf
    \State \textbf{Return} $\tilde\bx$, $\mathscr{I}^{\textrm{out}} = \{\mT^*, i+1, h\}$
\ElsIf{$\mT^k \neq \mathbb{R}^n$} 
    \State $h \gets h + 1$\label{step2:h_increase} \label{step:sub2_hupdate}
    \State $\mT^{k+1} \gets \text{S-UPDATE}(\mT^k)$
    \If {$h = h^* + 1$}
        \State $h \gets 0$ \label{step:sub2_hupdate_end}
    \EndIf
        \State $k \gets k+1$
    \State Return to Step \ref{step:REMOVEfunction}
\Else
    \State \textbf{Return} $\bhx$, $\mathscr{I}^{\textrm{in}}$
\EndIf
\end{algorithmic}
\end{subroutine}

We extend the previous cut generation subroutine to generate lower dimensional trust regions by randomly selecting $\mS \subset \mN$, as shown in Subroutine \ref{subroutine:CP-TRdimred}. At a high level, the main difference between Subroutine \ref{subroutine:CP-TR} and Subroutine \ref{subroutine:CP-TRdimred} is that when the trust region size ($p$) is increased, Subroutine \ref{subroutine:CP-TRdimred} first computes a series of low dimensional trust regions of the increased size, exhausts the corresponding subregions of violated forward-feasible points, and then considers the full dimensional trust region. More specifically, Subroutine \ref{subroutine:CP-TRdimred} considers the full dimensional trust region $\mT(\bhx, p)$ only after $h^*$ low dimensional trust regions of size $p$ have \emph{consecutively} failed to generate a cut. The index $h$ tracks the number of consecutive ``empty" low dimensional trust regions (step \ref{step2:h_increase}), and is reset to zero when a cut is successfully generated within a trust region (step \ref{step2:h_reset}). 

The trust region updating function is modified accordingly to reflect this process, and is slightly more involved than the original UPDATE function in equation \eqref{eq:update-centered}. To generate a low dimensional trust region, we randomly sample a set $\mS \subset \mN$. Let $s$ denote a value less than $|\mN|$, and let RAND($s$) denote a function that generates a random subset of $\mN$ of size $s$. The value $s$, which represents the cardinality of the set $\mS$, can be either a fixed, predetermined value or a function of the trust region size, i.e., $s(p)$. For now, we assume the latter, and will discuss further in the following paragraph. The new stochastic (indicated by S-) update function becomes:
\begin{equation}\label{eq:s-update}
\text{S-UPDATE}(\mT_{\mS}(\bhx,p)) = \begin{cases}
\mT(\bhx, p), & \quad \text{if } h = h^*,\\
\mT_{\text{RAND}(s(p))}(\bhx, p), & \quad \text{if } h < h^*,\\
\mT_{\text{RAND}(s(\delta p))}(\bhx, \delta p), & \quad \text{if } h = h^*+1.\\
\end{cases}
\end{equation}
where $\delta$, the trust region growth rate, is the same as in the original UPDATE function. When $h<h^*$, a new low dimensional trust region is generated. When $h=h^*$, the full dimensional trust region is generated instead. Finally, when the full dimensional trust region fails to identify any violated forward-feasible point, i.e., $h = h^*+1$, the S-UPDATE function generates a new low dimensional trust region of an increased size.  

We conclude by discussing the choice of the function $s(p)$ in equation \eqref{eq:s-update}, which defines the cardinality of the set $\mS$ to be sampled given that the trust region is of size $p$. While $s(p)$ can be defined simply to return a fixed value (less than $|\mN|$) for any input $p$, we define $s(p)$ to help offset the reduction in computational efficiency associated with larger trust regions. In particular, we define $s(p)$ to be a decreasing function of the trust region size $p$, i.e., we increasingly reduce the dimensionality of the trust region as its size grows. Our specific choice of the function $s(p)$ is:
\begin{equation}\label{eq:sampling_cardinality}
s(p) = \max\{\lfloor (1 - \kappa(p-1))n \rfloor, \lfloor qn \rfloor\}.  
\end{equation}
In equation~\eqref{eq:sampling_cardinality}, the parameter $\kappa \in (0,1)$ denotes the dimensionality reduction rate: when the size of the trust region increases by 1 unit, the dimension of the trust region decreases by $\kappa$. The parameter $q \in (0,1)$ specifies an upper bound on the dimensionality reduction, so that the dimension of the trust region will be no less than $qn$. 

\section{Extension to Multi-Point Inverse MILO}\label{sec:contemporary_models}

In this section, we present an extension of our base inverse MILO model~\eqref{model:IOM-Bilevel1}. While the literature surrounding inverse MILO problems has almost exclusively focused on model~\eqref{model:IOM-Bilevel1}, the literature on inverse convex optimization problems has expanded \iz{in the direction of data-driven parameter estimation}, where multiple solutions can be used as input. In light of this, we propose a \emph{multi-point} inverse MILO model, which takes multiple input solutions from potentially different forward-feasible regions to generate a cost vector. We show how our trust region-based cutting plane algorithms can be directly extended to solve this model.

\subsection{Problem description}

Let $\mD$ be the index set of $D$ data points, $\hat\bx_1, \ldots, \hat\bx_D$, which are feasible for their respective feasible regions $\mX_1, \ldots, \mX_D$. 
The following model is a natural extension of model~\eqref{model:IOM-Bilevel1} to the multi-point case:\iz{
\begin{subequations}\label{model:multi-point_IO}
\begin{align}
\underset{\bc \in \mP}{\text{minimize}} \quad & \norm{\bc-\bc^0}_1\\
\text{subject to} \quad & \bc \in \mC(\bhx_d, \mX_d), \quad \forall d \in \mD.\label{cons:invfeas_multi-point_a}
\end{align}
\end{subequations}
}
\iz{The data-driven inverse optimization literature has also considered ``noisy" data cases for continuous forward problems, where it is not possible to compute a single vector $\bc$ that is inverse-feasible for each and every decision input. In the light of this, we also propose the following extension of model \eqref{model:multi-point_IO}, where we introduce an inverse-feasible cost vector $\bbc_d$ for each data point $\bhx_d \in \mX_d$:}
\iz{
\begin{subequations}\label{model:multi-point_IO_noisy}
\begin{align}
\underset{\bc \in \mP, \{\bbc_d\}_{d \in \mD}}{\text{minimize}} \quad & \norm{\bc-\bc^0}_1 + \lambda \sum_{d \in \mD} \norm{\bbc_d-\bc}_1\\
\text{subject to} \quad & \bbc_d \in \mC(\bhx_d, \mX_d), \quad \forall d \in \mD \label{cons:invfeas_multi-point}.
\end{align}
\end{subequations}}
\iz{The regularization term $\lambda$ represents a trade-off between making the cost vector $\bc$ close to the prior $\bc^0$ and to the set of inverse-feasible cost vectors $\{\bbc_d\}_{d \in \mD}$.}

\subsection{Extended multi-point cutting plane algorithm}

We now show how a natural extension of Algorithm \ref{algo:CP-TR} can be used to solve models \iz{\eqref{model:multi-point_IO} and \eqref{model:multi-point_IO_noisy}}. The key observation is that each feasible region $\mX_d$ in \iz{\eqref{cons:invfeas_multi-point_a} and \eqref{cons:invfeas_multi-point}} can be replaced with a forward-feasible generator set $\mG_d \subseteq \mX_d$. So, the previous trust region ideas are applicable here.  The cutting plane algorithm \iz{for solving model \eqref{model:multi-point_IO_noisy}} is provided in Algorithm \ref{algo:CP-TRmulti}. \iz{The algorithm can be easily simplified to solve model \eqref{model:multi-point_IO} by ignoring all appearances of the set $\{\bbc_d\}_{d \in \mD}$ and replacing all appearances of terms $\bbc_d$ with $\bc$.}

\begin{algorithm}[h]
\caption{Cutting plane algorithm for multi-point inverse MILO}
\label{algo:CP-TRmulti}
\vspace{0.3cm}
\textbf{Input:} An inverse MILO problem instance $(\{\hat\bx_d, \mX_d \}_{d \in \mD})$, initial trust region size $p^0$\\
\textbf{Output:} An inverse-optimal solution $(\bc^{*}, \{\bbc^*_d\}_{d \in \mD})$
\begin{algorithmic}[1]
\State Initialize $i = 0, \mathscr{I}^0_d = (\mT(\bhx_k,p^0),i), \tilde\mX^i_d = \emptyset$,  count$_e$ = 0, count$_v$ = 0.
\While{count$_e$ $< |\mD|$}
\State count$_e$ $\gets$ 0
\State count$_v$ $\gets$ 0
 
\State Solve $\textbf{MP}(\{\bhx_d,\tilde\mX^i_d\}_{d \in \mD})$, let $\tilde\bc^{i}, \{\tilde\bc^i_d\}_{d \in \mD}$ be its optimal solution 
\While{(count$_v$ $< V^*$) and (count$_e$ + count$_v$ $< |\mD|$)} \label{step:multipoint_depthtolerance}
\For{$d = 1, \ldots ,D$} \label{step:multipoint_sort}
\State Solve \textbf{SUBROUTINE}($\tilde\bc^i_d, \bhx_d, \mX_d, \mathscr{I}^i_d)$, let $\tilde\bx^{i+1}_d$ and $\mathscr{I}^{i+1}_d$ be its output
\If {$\tilde\bx^{i+1}_d = \bhx_d$}
    \State count$_e$ $\gets$ count$_e$ + 1
    \State $\tilde\mX^{i+1}_d \gets \tilde\mX^{i}_d$
\Else
    \State count$_v$ $\gets$ count$_v$ + 1
    \State $\tilde\mX^{i+1}_d \gets \tilde\mX^{i}_d \cup \{\tilde\bx^{i+1}_d\}$
\EndIf
\EndFor
\EndWhile
\State $i \gets i+1$
\EndWhile \\
\Return $\bc^* = \tilde\bc^i, \ \bbc^*_d = \tilde\bc^i_d \ \forall d \in \mD$ 
\end{algorithmic}
\end{algorithm}

The master problem defines a relaxation of model \eqref{model:multi-point_IO_noisy}, where $\tilde\mX_d \subseteq \mX_d$ is a finite set for each $d$:
\iz{
\begin{subequations}\label{model:MP_multi-point_IO}
\begin{align}
\ \textbf{MP}(\bc^0,\{\bhx_d, \mX_d\}_{d \in \mD}) : \ \underset{\bc \in \mP, \{\bbc_d\}_{d \in \mD}}{\text{minimize}} \quad & \sum_{d\in\mD} \norm{\bc - \bc^0}_1 + \lambda \sum_{d \in \mD} \norm{\bbc_d-\bc}_1\\
\text{subject to} \quad & \bbc_d^\top(\bhx_d - \bx_d) \leq 0, \quad \forall \; \bx_d \in \tilde\mX_d, \ d \in \mD.
\end{align}
\end{subequations}}
With multiple data points, there is a cut generation problem associated with each $d \in \mD$. Thus, there is a new issue that must be considered, specifically that of a breadth-first versus depth-first search: whether the master problem is called after cuts for each $d$ are generated versus after the first cut for a single $d$. Let $V^* \in \{1, \ldots, D\}$ be a parameter that denotes the minimum number of violated forward-feasible points, if any, that must be found until the master problem is called. The larger $V^*$ is, the more the algorithm behaves ``breadth-first''. In the extreme case where $V^* = D$, the master problem is called only when the master problem candidate solution is considered for every $d$. This may reduce the number of master problem calls, but at the expense of potentially solving many unnecessary cut generation problems that cannot identify violated cuts. For example, this may happen when a generator set has been identified for several forward-feasible regions already but not all of them. In the other extreme where $V^* = 1$, the master problem is called as soon as a violated forward-feasible point is found for any $d \in \mD$. In the algorithm, count$_v$ (count$_e$) denotes the number of forward-feasible regions for which a violated point has (has not) been found. The algorithm terminates when count$_e = D$, i.e., there does not exist a single violated forward-feasible point in any forward-feasible region.


\section{Experimental Design}\label{sec:numerical_setup}

In this section, we discuss the setup of our computational experiments, which include the generation of test instances, the different subroutine variants and enhancements considered, and the parameterization of each subroutine. Our computational experiments are focused on inverse MILO problems of the form presented in model \eqref{model:IOM-Bilevel1} \iz{with $\mP = \mathbb{R}^n$}. 
\iz{We made this modeling choice to focus the presentation and discussion on the key insights gained from applying our new cutting plane algorithms over various problem structures and forward-feasible regions. This model is also studied in all previous computational studies of inverse MILO \citep{wang2009cutting, schaefer2009inverse, duan2011heuristic}.} To generate comprehensive insights into the performance of our algorithms, we apply them to a diverse set of problem instances, much larger than any set considered in the literature to date.

\subsection{Test bank generation}

To generate a comprehensive set of inverse MILO problems over many different problem structures, we draw from the MIPLIB 2017 benchmark library \citep{gleixner2019miplib}, which includes 240 MIPLIB problems. We use a subset of these problems to generate a test bank of inverse MILO problems as follows. We first consider all MIPLIB problems with fewer than 12,000 variables and constraints. This reduces our bank to 125 total MIPLIB problems. For each problem, we attempt to generate a feasible solution by solving the problem with a randomly generated cost vector. If the problem is unbounded or cannot be solved within 10 minutes using Gurobi 8.1.0, we try again with a different randomly generated cost vector. If this procedure cannot generate three feasible solutions within ten attempts, the MIPLIB problem is dropped from consideration. Otherwise, three feasible points are generated for each MIPLIB problem. This approach results in 73 MIPLIB problems being included in our test bank, each with three feasible points, totalling 219 inverse MILO problem instances. 
Each instance is constructed by using one of the feasible solutions as $\bhx$, the corresponding original cost vector of the MIPLIB problem (i.e., not the random one used to generate $\bhx$) as $\bc^0$ and the corresponding set of MIPLIB problem constraints as $\mX$. The 3 instances constructed per MIPLIB problem are labeled using the MIPLIB name and a suffix of t1, t2, and t3 (see Table \ref{TableEC:all_results}).


Of the 73 included MIPLIB problems, 9 have only binary variables, 42 have binary and continuous variables, 14 have binary and integer variables, and 8 have binary, integer and continuous variables. The number of variables in these problems ranges from 34 to 11,700 with an average of 3,185. The number of constraints ranges from 4 to 10,900 with an average of 2,315. These problems represent applications such as scheduling, production, shipment, assignment, set covering and bin packing. 


\subsection{Subroutine variants}
To solve these inverse MILO instances, we implement and compare five different cutting plane algorithms, described below:

\begin{itemize}

\item \textbf{CP.} This is the classical cutting plane approach proposed in \cite{wang2009cutting}, which serves as the benchmark algorithm. 

\item \textbf{CP-ES.} This is a simple extension of the CP algorithm in which we embed the early-stop heuristic presented in Section~\ref{subsec:early-stop}. We choose a time threshold of 5 seconds, at which time the most violated forward-feasible solution is returned ($\tau = 5$). If no violated points have been found by 5 seconds, the subroutine returns the first one found after the 5-second threshold has been exceeded.

\item \textbf{CPTR.} This is the cutting plane algorithm with trust regions, defined in Subroutine \ref{subroutine:CP-TR}. The subroutine initializes with a trust region of size one ($p^0 = 1$ in equation \eqref{eq:normed_TR}). The trust region doubles in size each time the trust region becomes redundant ($\delta = 2$ in equation \eqref{eq:update-centered}). The trust region is removed every 10th cut that is generated and when no violated forward-feasible points can be found within the trust region ($i^* = 10$ and $k^* = 2$ in equation \eqref{eq:eval_function}).

\item \textbf{CPTR-ES.} This is an extension of the CPTR algorithm in which we embed the early-stop heuristic with a threshold of $5$ seconds ($\tau = 5$), similar to the CP-ES algorithm.

\item \textbf{CPTR-ES-DR.} This considers the CPTR algorithm with both the early-stop heuristic and the stochastic dimensionality reduction enhancement. We implement a dimensionality reduction rate of 3\% ($\kappa = 0.03$ in equation \eqref{eq:sampling_cardinality}). The dimension of the trust region will never be lower than 80\% of $n$ ($q$ = 0.8). Finally, we revert to the full dimensional trust region when ten consecutive low-dimensional trust regions are found to have no violated forward-feasible points ($h^* = 10$). 
\end{itemize}

All experiments were coded in Python 3.7 and optimization problems are solved using Gurobi 8.1.0 parameterized with a single thread. The experiments were conducted on a Intel Core i7-4790 processor at 3.60 GHz on a Windows 8.1 Pro. A time constraint of one hour was set for all instances. 


\section{Numerical Results}\label{sec:numerical_ex}

\begin{figure}[h]
\centering
\includegraphics[width = 0.65\textwidth]{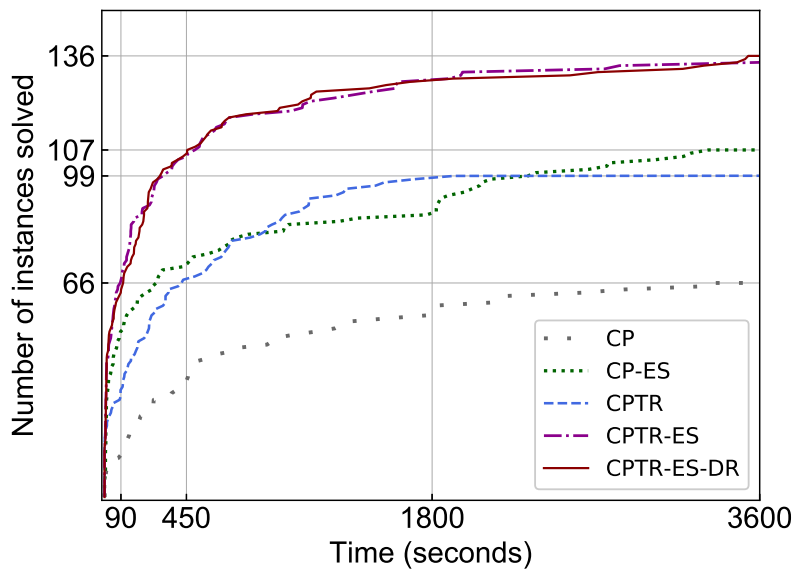}
\caption{Performance profiles of different cutting-plane algorithms\iz{; 66, 99, 107, 134 and 136 instances, out of 215 total instances, were solved by the CP, CP-ES, CPTR, CPTR-ES and CPTR-ES-DR algorithms respectively over a one-hour time window}.}\label{fig:general_results}
\end{figure}

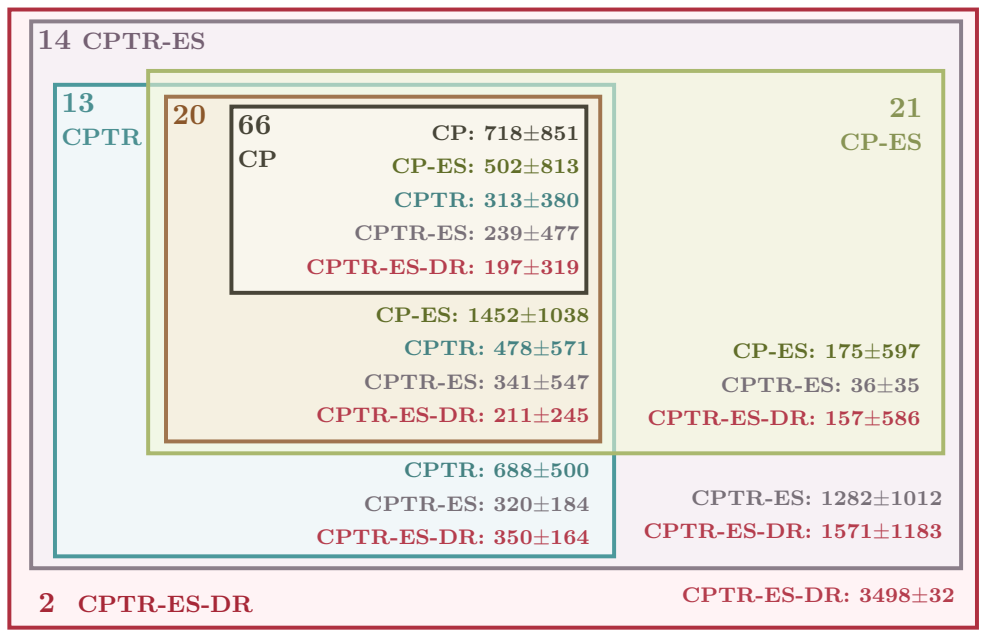
\begin{figure}[h]
\linespread{1}
\centering
\begin{tikzpicture}[x=0.75pt,y=0.75pt,yscale=-1,xscale=1]

\draw  [color={rgb, 255:red, 175; green, 48; blue, 64 }  ,draw opacity=1 ][fill={rgb, 255:red, 255; green, 239; blue, 241 }  ,fill opacity=0.7 ][line width=1.5]  (142,23) -- (630,23) -- (630,335) -- (142,335) -- cycle ;
\draw  [color={rgb, 255:red, 139; green, 127; blue, 138 }  ,draw opacity=1 ][fill={rgb, 255:red, 246; green, 241; blue, 246 }  ,fill opacity=0.8 ][line width=1.5]  (153,29) -- (622,29) -- (622,305) -- (153,305) -- cycle ;
\draw  [color={rgb, 255:red, 76; green, 153; blue, 156 }  ,draw opacity=1 ][fill={rgb, 255:red, 237; green, 252; blue, 253 }  ,fill opacity=0.55 ][line width=1.5]  (165,61) -- (447,61) -- (447,299) -- (165,299) -- cycle ;
\draw  [color={rgb, 255:red, 170; green, 184; blue, 111 }  ,draw opacity=1 ][fill={rgb, 255:red, 242; green, 248; blue, 215 }  ,fill opacity=0.55 ][line width=1.5]  (212,54) -- (613,54) -- (613,247) -- (212,247) -- cycle ;
\draw  [color={rgb, 255:red, 139; green, 87; blue, 42 }  ,draw opacity=0.8 ][fill={rgb, 255:red, 250; green, 235; blue, 222 }  ,fill opacity=0.55 ][line width=1.5]  (221,67) -- (440,67) -- (440,241) -- (221,241) -- cycle ;
\draw  [color={rgb, 255:red, 72; green, 69; blue, 56 }  ,draw opacity=1 ][fill={rgb, 255:red, 250; green, 248; blue, 241 }  ,fill opacity=0.8 ][line width=1.5]  (254,72) -- (433,72) -- (433,166) -- (254,166) -- cycle ;

\draw (290.5,80.17) node [anchor=north west][inner sep=0.75pt]   [align=right] {{\scriptsize \textcolor[rgb]{0.28,0.27,0.22}{\textbf{CP: 718$\pm$851 }}}\\{\scriptsize \textcolor[rgb]{0.38,0.43,0.16}{\textbf{CP-ES: 502$\pm$813 }}}\\{\scriptsize \textcolor[rgb]{0.27,0.51,0.51}{\textbf{CPTR: 313$\pm$380 }}}\\{\scriptsize \textcolor[rgb]{0.47,0.44,0.47}{\textbf{CPTR-ES: 239$\pm$477 }}}\\{\scriptsize \textcolor[rgb]{0.71,0.24,0.3}{\textbf{CPTR-ES-DR: 197$\pm$319 }}}};
\draw (256,75) node [anchor=north west][inner sep=0.75pt]  [color={rgb, 255:red, 72; green, 69; blue, 56 }  ,opacity=1 ] [align=left] {\textbf{66}\\{\footnotesize \textbf{CP}}};
\draw (167,64) node [anchor=north west][inner sep=0.75pt]  [color={rgb, 255:red, 76; green, 144; blue, 146 }  ,opacity=1 ] [align=left] {\textbf{13}\\{\footnotesize \textbf{CPTR}}};
\draw (559.5,66.5) node [anchor=north west][inner sep=0.75pt]  [color={rgb, 255:red, 136; green, 148; blue, 85 }  ,opacity=1 ] [align=right] {\textbf{21}\\{\footnotesize \textbf{CP-ES}}};
\draw (155,32) node [anchor=north west][inner sep=0.75pt]  [color={rgb, 255:red, 117; green, 96; blue, 115 }  ,opacity=1 ] [align=left] {\textbf{\textcolor[rgb]{0.46,0.38,0.45}{14 }{\footnotesize \textcolor[rgb]{0.46,0.38,0.45}{CPTR-ES}}}};
\draw (295.5,250.5) node [anchor=north west][inner sep=0.75pt]   [align=right] {{\scriptsize \textcolor[rgb]{0.27,0.51,0.51}{\textbf{CPTR: 688$\pm$500 }}}\\{\scriptsize \textcolor[rgb]{0.47,0.44,0.47}{\textbf{CPTR-ES: 320$\pm$184 }}}\\{\scriptsize \textcolor[rgb]{0.71,0.24,0.3}{\textbf{CPTR-ES-DR: 350$\pm$164 }}}};
\draw (462.5,190.5) node [anchor=north west][inner sep=0.75pt]   [align=right] {{\scriptsize \textcolor[rgb]{0.38,0.43,0.16}{\textbf{CP-ES: 175$\pm$597}}}\\{\scriptsize \textcolor[rgb]{0.47,0.44,0.47}{\textbf{CPTR-ES: 36$\pm$35}}}\\{\scriptsize \textcolor[rgb]{0.71,0.24,0.3}{\textbf{CPTR-ES-DR: 157$\pm$586}}}};
\draw (460.5,264.5) node [anchor=north west][inner sep=0.75pt]   [align=right] {{\scriptsize \textcolor[rgb]{0.47,0.44,0.47}{\textbf{CPTR-ES: 1282$\pm$1012 }}}\\{\scriptsize \textcolor[rgb]{0.71,0.24,0.3}{\textbf{CPTR-ES-DR: 1571$\pm$1183 }}}};
\draw (295.5,172) node [anchor=north west][inner sep=0.75pt]   [align=right] {{\scriptsize \textcolor[rgb]{0.38,0.43,0.16}{\textbf{CP-ES: 1452$\pm$1038 }}}\\{\scriptsize \textcolor[rgb]{0.27,0.51,0.51}{\textbf{CPTR: 478$\pm$571 }}}\\{\scriptsize \textcolor[rgb]{0.47,0.44,0.47}{\textbf{CPTR-ES: 341$\pm$547 }}}\\{\scriptsize \textcolor[rgb]{0.71,0.24,0.3}{\textbf{CPTR-ES-DR: 211$\pm$245 }}}};
\draw (223,70) node [anchor=north west][inner sep=0.75pt]  [color={rgb, 255:red, 139; green, 87; blue, 42 }  ,opacity=1 ] [align=left] {\textbf{20}};
\draw (502.5,82.5) node [anchor=north west][inner sep=0.75pt]  [font=\small,color={rgb, 255:red, 136; green, 148; blue, 85 }  ,opacity=1 ] [align=left] {};
\draw (475.5,313.5) node [anchor=north west][inner sep=0.75pt]   [align=right] {{\scriptsize \textcolor[rgb]{0.71,0.24,0.3}{\textbf{ CPTR-ES-DR: 3498$\pm$32 }}}};
\draw (155.5,316.5) node [anchor=north west][inner sep=0.75pt]  [color={rgb, 255:red, 167; green, 39; blue, 55 }  ,opacity=1 ] [align=left] {\textbf{2 \ {\footnotesize CPTR-ES-DR}}};

\end{tikzpicture}
\caption{Summary of computational results over solved instances. Single numbers denote the number of solved instances in the difference sets. Solution times in seconds are given in the form of average $\pm$ standard deviation.}\label{fig:performance_summary}
\end{figure}

In this section, we present the computational results using the algorithms outlined in Section \ref{sec:numerical_setup}. We begin with a high-level overview of the results, summarized by Figures \ref{fig:general_results} and \ref{fig:performance_summary}. We then provide an in-depth study of the advantages and potential trade-offs that exist when using trust regions and the additional enhancements (Section \ref{subsec:results_analysis}). In particular, we compare cut generation speed and cut strength across the various algorithms. We also conduct a sensitivity analysis, examining how the results change when several baseline parameters are modified (Section \ref{subsec:baseline_values}). We conclude with a brief summary of the main takeaways from our numerical results (Section \ref{subsec:computational_summary}).

From the performance profiles in Figure \ref{fig:general_results}, we observe that CPTR-ES and CPTR-ES-DR solve significantly more instances than CPTR or CP-ES, which in turn solve many more than CP, the baseline algorithm. The two best algorithms are also significantly faster than the others. For instance, they can each solve the same number of instances as CP  (66 instances) and CP-ES (107 instances) in only 2.5\% and 12\% of the required time, respectively. The main takeaway from these high-level results is that either trust regions or the early-stop heuristic can lead to dramatic improvements in the solution time, compared to CP. When trust regions and early-stopping are used together, the improvements are even more significant.

\linespread{1}

Figure~\ref{fig:performance_summary} presents a more granular breakdown of the solved instances. Notably, the results exhibit a distinct \emph{nested} structure in which each additional enhancement can solve strictly more instances than without the enhancement. For example, all 66 instances solved by CP are also solved by CP-TR and CP-ES, which solve an additional 33 and 41 instances, respectively. CPTR-ES solves 134 instances, including all 120 instances solved by at least one of CP-ES or CPTR. Finally, CPTR-ES-DR solves 2 more instances beyond those solved by CPTR-ES. Within each group of solved instances, CPTR-ES and CPTR-ES-DR have the fastest solution times. An example of algorithm performance over a particular instance is shown in Section \ref{EC:example_of_instance} in the Electronic Companion.

Before we proceed to a more detailed discussion of these results, we acknowledge that 83 of the instances were not solved by any of the algorithms, including our best-performing ones. Ultimately, the need to solve many MILO problems is an inherent limitation of all these cutting plane algorithms. The MILO problems are used to generate cuts and validate optimality of a candidate cost vector. As these are MIPLIB problems, many of the MILO problems that need to be solved are time consuming relative to our chosen time limit. While trust regions can make certain MILO problems easier to solve, the algorithms eventually need to solve the full MILO problem to verify optimality. And we observe that many of the unsolved instances are precisely the ones in which MILO solution times increase dramatically when larger trust regions are considered.

\subsection{Analysis and discussion of results}\label{subsec:results_analysis}

In this section, we analyze how the addition of trust regions and other enhancements lead to the observed performance results. The analysis is presented in the following order: (i) examining the effects of the early-stop heuristic on CP, (ii) analyzing the gains from adding trust regions to both CP and CP-ES, (iii) examining the effects of adding dimensionality reduction on CPTR-ES.  

\begin{figure}[h]
     \centering
     \includegraphics[width = 0.87\textwidth]{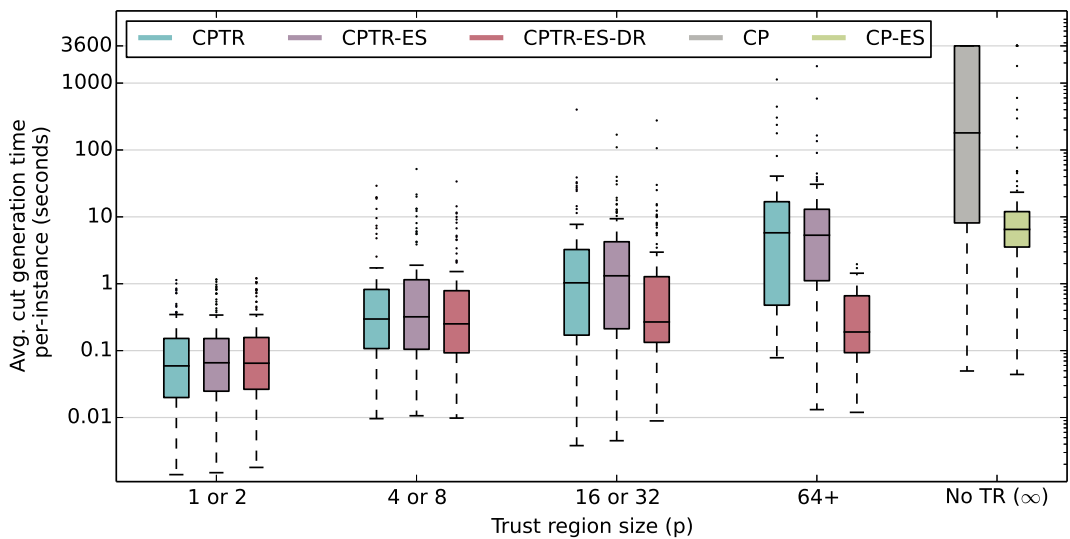}
     \caption{The average per-instance cut generation times within different trust region sizes and with enhancements. Trust region sizes correspond to $p$ in equation~\eqref{eq:normed_TR} and \eqref{eq:TR_with_DR}. Each data point is the average cut generation time within a trust region size for a particular instance. \iz{$ p = 64+$ denotes all finite trust region sizes of $p \geq 64$, whereas} No TR denotes no trust region (i.e., $p = \infty$).}
     \label{fig:avg_CGT}
 \end{figure}

The two main factors that determine overall solution times are cut generation speed and strength of cuts. Solution times are reduced when both factors are improved together, or when the improvement in one eclipses any potential loss in the other. To facilitate the following discussion, we use Figure \ref{fig:avg_CGT} to highlight the average cut generation speed over different trust region sizes and enhancements. In each of the three following subsubsections, we draw from Figure \ref{fig:avg_CGT} and additional information about the strength of cuts (tailored for each section) to analyze the results. In this subsection, all values appearing in figures are plotted on a log scale.

\subsubsection{Early-stop heuristic.} 
Cut generation times within CP can be very high, as evident from Figure \ref{fig:avg_CGT}. The median of the average per-instance cut generation times exceeds 200 seconds, and CP fails to generate any cut within the time limit in over 25\% of the instances. When the early-stop heuristic is applied, the median is reduced to 7 seconds, close to the early-stop $\tau$ value of $5$ seconds. CP-ES is thus capable of generating many more cuts, as described next.

\begin{figure}[h]
    \centering
    \includegraphics[width = 8.1cm]{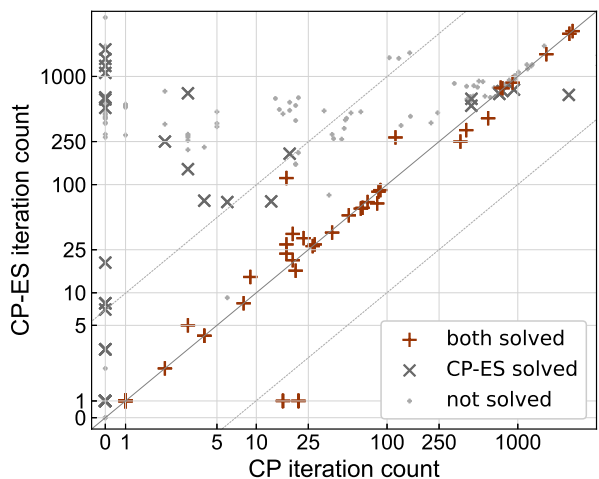}
    \caption{Iteration count and instance solvability when early-stop heuristic is added to CP.}\label{fig:CPsolved}
\end{figure}

The iteration count for instances solved by CP-ES is shown in Figure \ref{fig:CPsolved}. First, the ability to generate cuts faster allows CP-ES to generate many more cuts, which enables it to solve many more instances. Second, we observe that compared to the instances solved by CP, the early-stop heuristic does not significantly increase the iteration count. Part of this reason is that in a number of instances, CP and CP-ES are equivalent because the early-stop feature is never used; many of the instances that can be solved by CP are naturally ones in which cuts can be computed quickly. Nonetheless, the results in general suggest that cuts generated by the early-stop heuristic with $\tau = 5$ are not significantly weaker. 

Note that lowering the value of $\tau$ will not necessarily result in faster cut generation. For example, there exist many instances in which the average cut generation time is higher than $\tau$ (as observed in Figure \ref{fig:avg_CGT}), i.e., it takes longer than $\tau$ seconds on average to find any violated forward-feasible point. Cut generation times can also vary significantly depending on the given cost vector. In particular, many instances with an average cut generation time less than $\tau$ may still have a large number of iterations in which cut generation times are higher than $\tau$. 

\subsubsection{Trust regions.}

\begin{figure}[h]
    \centering
    \hspace{-0.35cm}
    \begin{subfigure}[t]{0.49\textwidth}
    \includegraphics[width = 8.1cm]{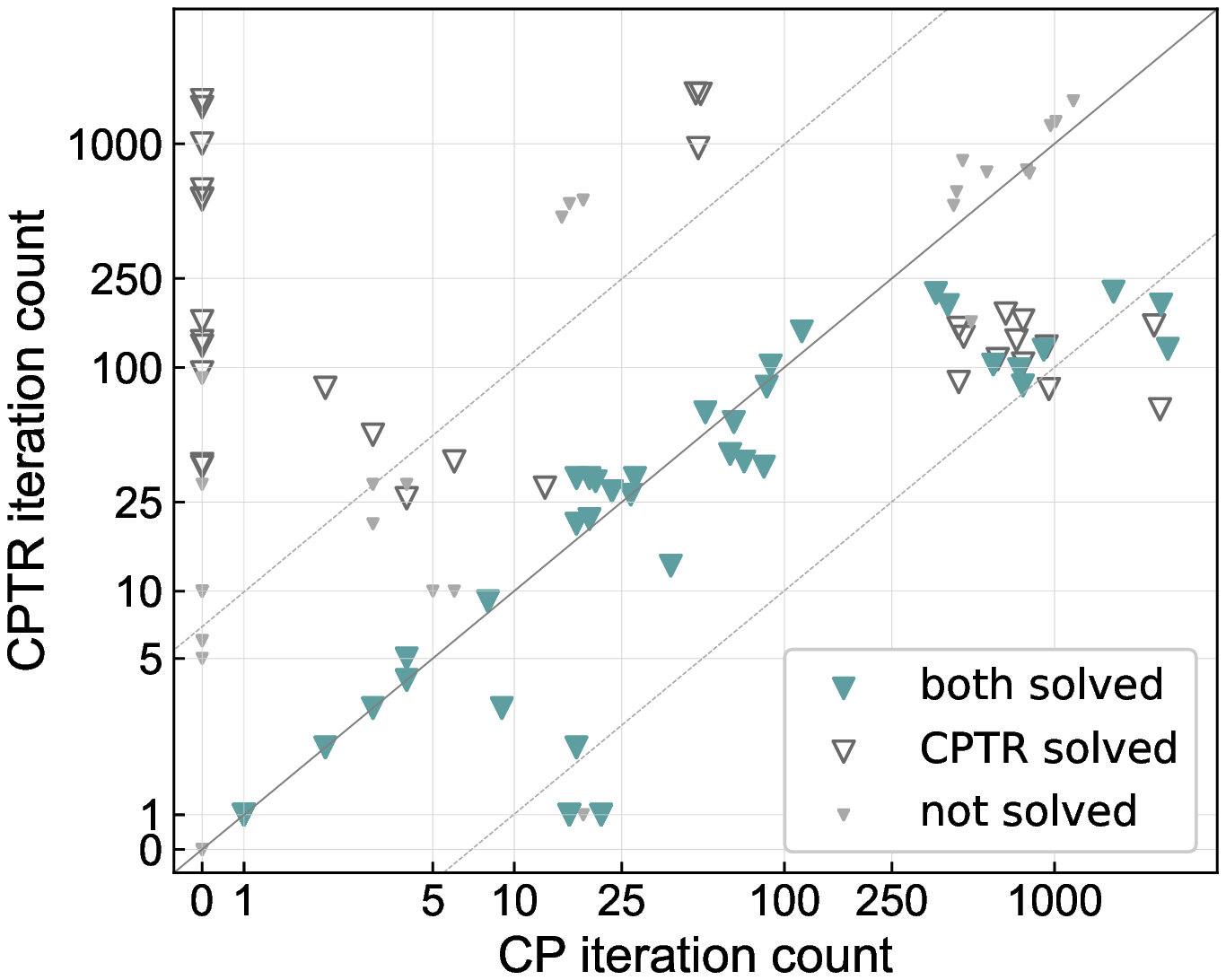}
    \caption{CP versus CPTR}\label{fig:TRvsCP}
    \end{subfigure}\ \
    \begin{subfigure}[t]{0.50\textwidth}
    \includegraphics[width = 8.2cm]{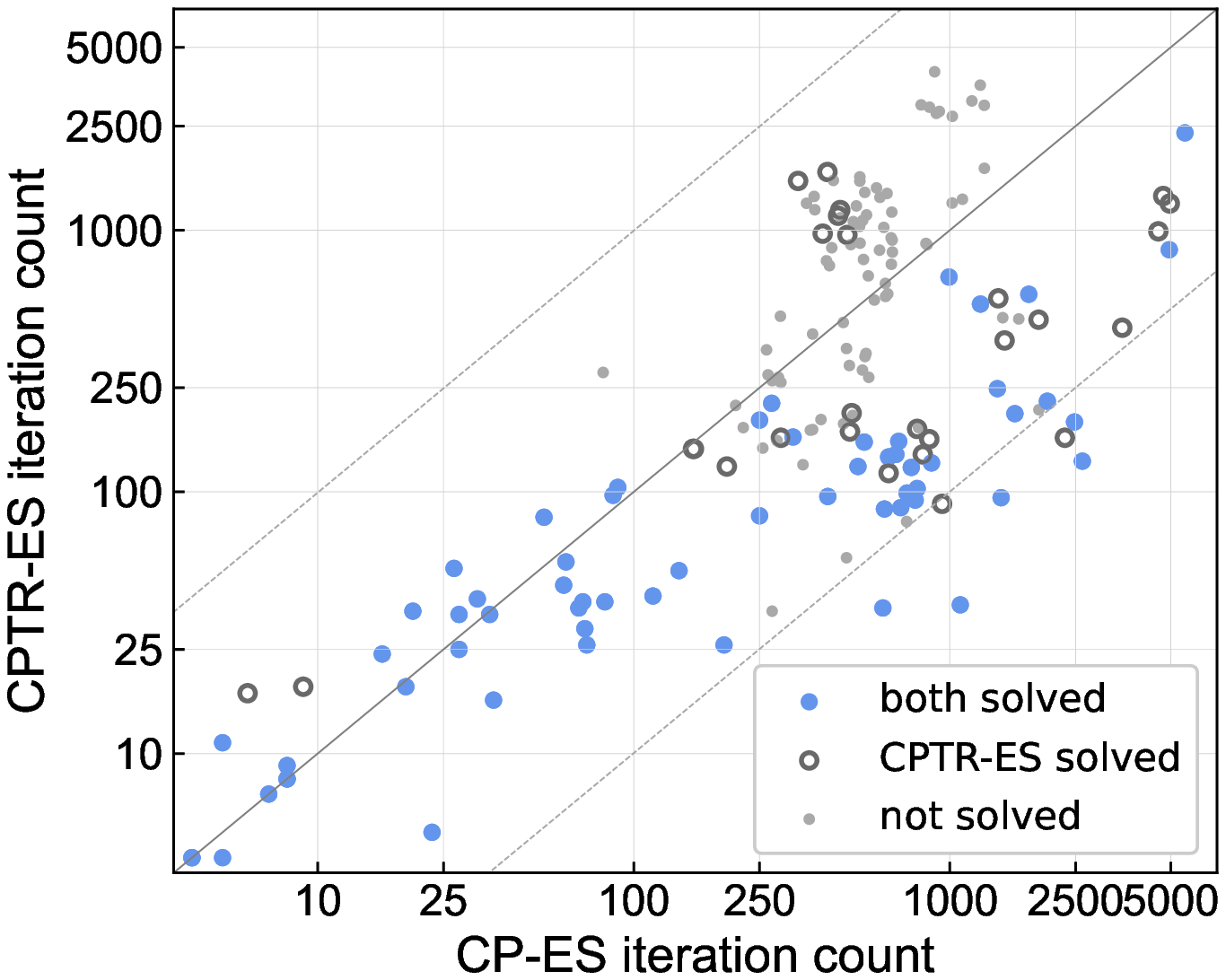}
    \caption{CP-ES versus CPTR-ES}\label{fig:TRESvsES}
    \end{subfigure}
    \vspace{0.2cm}
    \caption{Iteration count and instance solvability after adding trust regions to the CP and CP-ES algorithms.}\label{fig:iteration_count}
\end{figure}

Cut generation times can be substantially reduced using trust regions, as observed in Figure \ref{fig:avg_CGT}. For example, average cut generation time within trust regions of size $p \leq 8$ can be orders of magnitudes lower than without trust regions. Like the early-stop heuristic, this reduction in cut generation time allows more cuts to be generated when instances are unsolved by CP, and the improved cut generation speed leads to many more instances being solved, as seen in Figure \ref{fig:TRvsCP}. Of course, as trust region sizes increase, we observe that the cut generation times also tend to increase. 

For the instances solved by CP, CPTR solves those same instances in significantly fewer iterations. This observation highlights that the cuts computed by CPTR are stronger. For example, the instances that require over 250 cuts from CP require around an order of magnitude fewer cuts with CPTR. Comparing CPTR-ES to CP-ES, the reduction in iterations is similar, as seen in Figure \ref{fig:TRESvsES}. In general, the more cuts CP and CP-ES require to solve an instance, the greater the impact of trust regions on reducing iteration count.

The reduction in cut generation time and the increase in strength of cuts make trust regions very effective. Nonetheless, the magnitude of improvement in these two factors varies across instances. For example, reduction in cut generation time depends on how many cuts are computed within trust regions of smaller sizes. Secondly, a trust region of a particular size may lead to much faster cut generation in one instance than when the same trust region is applied to a different instance. To understand this phenomenon, it is worth noting that trust regions are, by definition, an additional set of constraints on the forward optimization problem. In some MIPLIB instances, the addition of trust regions of any size can ``simplify" the forward-feasible region and lead to faster cut generation, whereas in other instances large trust regions may make the problem even more difficult to solve.

\begin{table}[h]
\centering
\resizebox{!}{5.8cm}{%
\setlength{\tabcolsep}{12pt}
{\renewcommand{\arraystretch}{1.05}%
 \begin{tabular}{l r r r r  r r r}
 \toprule
& \multicolumn{2}{c}{CP-ES } & \multicolumn{2}{c}{CPTR-ES} & \multicolumn{3}{c}{Change (\%)} \\ \cmidrule(lr){2-3}\cmidrule(lr){4-5}\cmidrule(lr){6-8}
Instance & num.it. & time (s) & num.it. & time (s) & num.it. & \hspace{-0.5cm}c.g.time & \hspace{-0.5cm} time \\[0.1cm]
\midrule
\multicolumn{8}{l}{\textit{Instances that benefit from reduced iteration count}} \\[0.2cm]
mzzv42z\_t3 & 945 & $>$3600 & 90 & 285 & -90 & -18 & -92 \\
mzzv11\_t2 & 732 & 3181 & 99 & 404 & -86 & -7 & -87 \\
drayage-25-23\_t2 & 875 & 2657 & 129 & 321 & -85 & -19 & -88 \\
csched008\_t1 & 23 & 260 & 5 & 64 & -78 & -2 & -75 \\
sp150x300d\_t2 & 2035 & 181 & 222 & 29 & -89 & 46 & -84 \\[0.2cm]
\midrule
\multicolumn{8}{l}{\textit{Instances that benefit from reduced cut generation time}} \\[0.2cm]
seymour1\_t2 & 396 & $>$3600 & 970 & 263 & 145 & -97 & -93 \\
seymour\_t1 & 443 & $>$3600 & 1136 & 449 & 156 & -95 & -88 \\
roi2alpha3n4\_t2 & 197 & $>$3600 & 125 & 524 & -37 & -77 & -85 \\
neos-3083819-nubu\_t3 & 250 & 1837 & 188 & 343 & -25 & -75 & -81 \\
ran14x18-disj-8\_t1 & 994 & 726 & 662 & 96 & -33 & -80 & -87 \\[0.2cm]
\midrule
\multicolumn{8}{l}{\textit{Instances that benefit from both factors}} \\[0.2cm]
neos-4954672-berkel\_t1 & 1909 & $>$3600 & 455 & 258 & -76 & -70 & -93 \\
50v-10\_t2 & 4567 & $>$3600 & 990 & 121 & -78 & -85 & -97 \\
drayage-25-23\_t1 & 861 & $>$3600 & 159 & 270 & -82 & -60 & -92 \\
drayage-100-23\_t1 & 788 & $>$3600 & 174 & 274 & -78 & -66 & -92 \\
csched008\_t3 & 193 & 2768 & 26 & 130 & -87 & -66 & -95 \\
neos5\_t3 & 81 & 2760 & 38 & 0 & -53 & -100 & -100 \\
glass-sc\_t3 & 640 & 2074 & 136 & 146 & -79 & -67 & -93 \\
ran14x18-disj-8\_t2 & 1777 & 1289 & 569 & 95 & -68 & -77 & -93 \\
neos-4338804-snowy\_t2 & 1079 & 652 & 37 & 8 & -97 & -65 & -99\\
mik-250-20-75-4\_t1 & 2626 & 624 & 131 & 7 & -95 & -78 & -99 \\
\bottomrule
\end{tabular}}}\\
\vspace{0.2cm}
\caption{A subset of instances which illustrate the two key benefits of trust regions. The number of iterations (num.it.) and total time in seconds are provided, as well as the percentage in change of CPTR-ES over CP-ES in the number of iterations, cut generation time (c.g.time), and total solution time.}
\label{table:benefits}
\end{table}

In Table \ref{table:benefits} we show the results of adding trust regions to CP-ES over a subset of instances (selected from Table \ref{TableEC:all_results}). The results are partitioned into three categories: ones that benefit mainly from improved cut strength, ones that benefit mainly from reduced cut generation time, and ones that benefit from both. While the results shown are for some of the better performing instances, they show performance over many different MIPLIB problem structures and paint a representative picture of the overall results. Specifically, in the majority of instances, average cut generation times are reduced, which alone is sufficient in reducing solution time. When coupled with improved cut strength, even larger reductions in solution time can be expected. In instances where average cut generation time is not reduced for the reasons described in the previous paragraph, reduced solution times can still be expected as a result of reduced iteration count. Finally, we note that the comparison is made over CP-ES, which already performs significantly better than CP, the baseline algorithm.


\subsubsection{Dimensionality reduction.} 
\begin{figure}[h]
    \centering
    \includegraphics[width = 8.4cm]{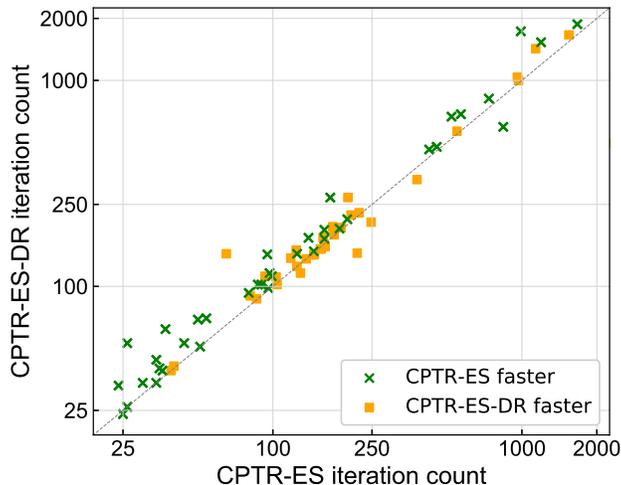}
    \caption{Comparison between the CPTR-ES and CPTR-ES-DR over instances that are solved by both.}\label{fig:CPTR-ESvsESS}
\end{figure}

Figure \ref{fig:CPTR-ESvsESS} shows the iteration count and overall performance over the 134 instances solved by both the CPTR-ES and CPTR-ES-DR algorithms. 

The first observation is that with very few exceptions, CPTR-ES-DR requires more iterations to solve each instance. This illustrates that considering lower dimensional trust regions generally leads to weaker cuts. The second observation illustrates the distinct trade-off that exists between cut generation speed and cut strength: when CPTR-ES-DR takes too many additional iterations, depicted by points being far to the left of the diagonal in Figure \ref{fig:CPTR-ESvsESS}, CPTR-ES generally performs better. When CPTR-ES and CPTR-ES-DR take a similar number of iterations, i.e., points close to the diagonal, CPTR-ES-DR performs better. While CPTR-ES-DR does lead to more iterations in general, if this increase is modest, the reduction in cut generation time will lead to better overall performance. On the other hand, if the increase in iteration count is too large, CPTR-ES will generally perform better. The parameters of the dimensionality reduction enhancement can be tuned to balance this trade-off, which we discuss in the next subsection. 
 
\subsection{Sensitivity analysis}\label{subsec:baseline_values}

In the previous section, we analyzed results for a fixed set of baseline parameters to illustrate that performance gains can be achieved over a wide variety of problem structures without tuning. In this section, we consider several modifications of the baseline parameter values and show that even minor algorithm tuning can lead to dramatic gains. We consider modifications to most parameter values, as shown in Table \ref{table:baseline_parameters}, and report the change in performance over the baseline values. For simplicity, we consider each parameter modification independently, holding all other parameters fixed to their baseline values. We choose a convenience sample of the ``\_t1''  instances solved by both the CPTR-ES and CPTR-ES-DR algorithms using baseline parameters, resulting in a total of 42 instances. 

\begin{table}[h]
\centering
\setlength{\tabcolsep}{10pt}
{\renewcommand{\arraystretch}{1.1}%
\begin{tabular}{ccccccccc}
 \toprule
  & \multicolumn{4}{c}{TR} & ES & \multicolumn{3}{c}{DR}\\
 \cmidrule(lr){2-5}\cmidrule(lr){6-6}\cmidrule(lr){7-9}
 Algorithm Parameters & $p^0$ & $\delta$ & $k^*$ & $i^*$ & $\tau (s)$ & $\kappa$ & $h^*$ & $q$ \\[0.1cm]
 \hline\\[-0.2cm]
& & 1.5 &  & 5 & 1 & 0.01 & 5 & 0.7\\
 \rowcolor{gray!20} Baseline Values & 1& 2 & 2 & 10 & 5 & 0.03 & 10 & 0.8 \\
& & 4 &  & 25 & 20 & 0.05 & 15 & 0.9\\
 \bottomrule
\end{tabular}}
\caption{List of considered algorithm parameters for the CPTR algorithm and the two enhancements.}\label{table:baseline_parameters}
\end{table}

\begin{figure}[h]
    \centering
    \begin{subfigure}[t]{0.48\textwidth}
    \includegraphics[width = 7.6cm]{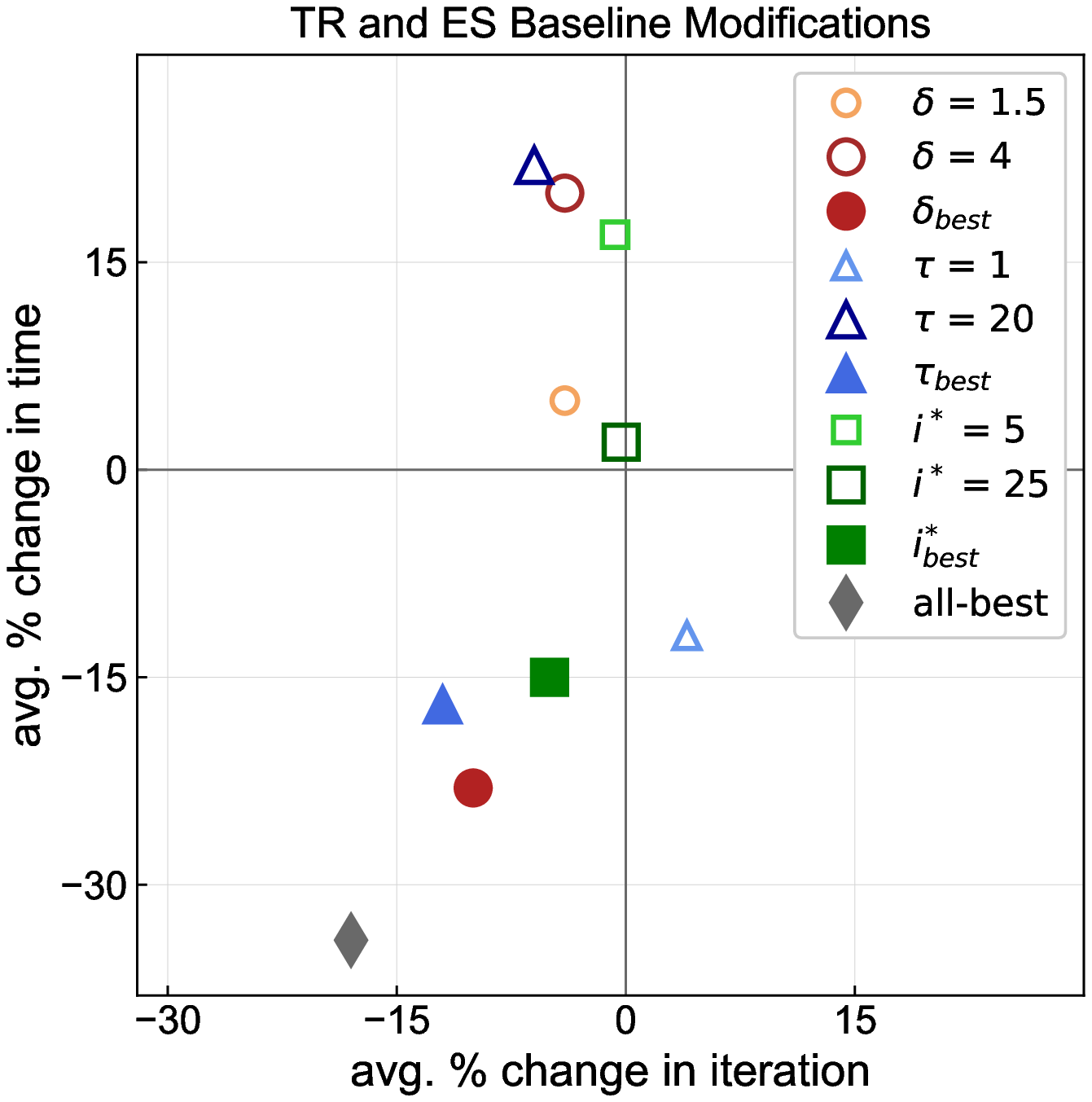}
    \caption{CPTR-ES}\label{fig:CPTR-ES-tuning}
    \end{subfigure}
    \begin{subfigure}[t]{0.48\textwidth}
    \includegraphics[width = 7.6cm]{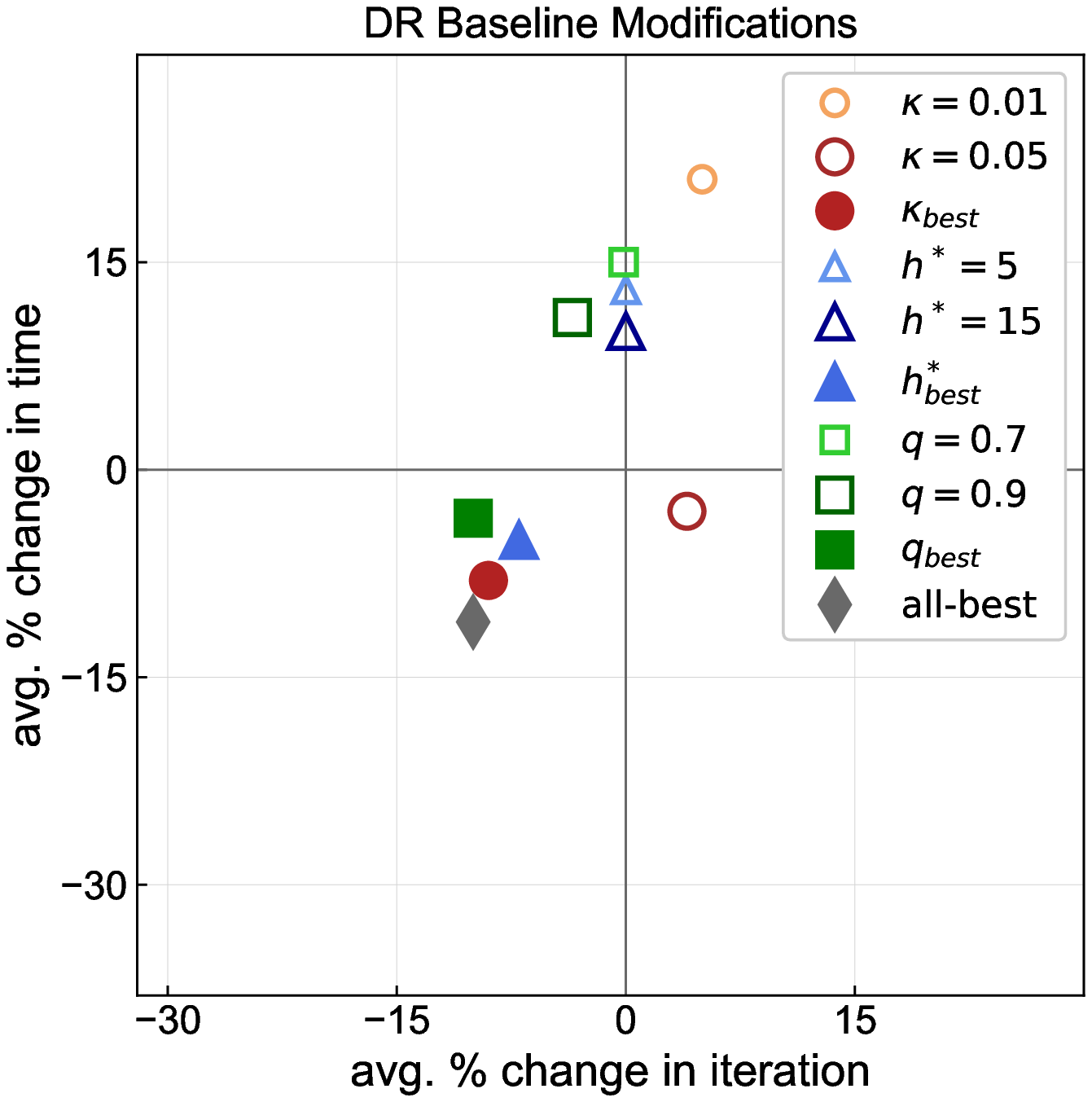}
    \caption{CPTR-ES-DR}\label{fig:CPTR-ES-DR-tuning}
    \end{subfigure}
    \vspace{0.2cm}
    \caption{Modification of parameter values for the CPTR-ES and CPTR-ES-DR. Baseline values are $\delta = 2, \tau = 5, i^* = 10$, $\tau = 5$, $\kappa = 0.03$, $h^* = 10$, $q = 0.8$. The ``best" subscript denotes the results when the best of the three (including baseline) parameter values is chosen for each instance, whereas the ``all-best" label corresponds to the case where the best of all parameter values is chosen. }\label{fig:sensitivity}
\end{figure}

The percentage change in total iteration count and solution time for CPTR-ES and CPTR-ES-DR, relative to the baseline parameter values, are shown in Figure \ref{fig:sensitivity}. 
The results highlight that, on average, the baseline parameter values perform quite well relative to the modified parameter values. However, we also show that it is possible to achieve significantly better performance if we could identify the best parameter values for each instance independently (markers denoted ``best''). For example, for CPTR-ES, if we could choose the best values for all three parameters for all instances, we could further reduce solution times by over 33\%. Similarly, for CPTR-ES-DR, it is possible to achieve an 11\% reduction in solution time with the best parameter values for each instance. While it is not possible to know the best parameters for each instance in advance, the results demonstrate that our default parameters tend to work well and that some tuning could lead to further improvements.

Comparing Figures \ref{fig:CPTR-ES-tuning} and \ref{fig:CPTR-ES-DR-tuning}, it appears that solution time is more sensitive to the parameters for trust regions and the early-stop heuristic than the parameters for the dimensionality reduction enhancement. This result is intuitive since trust regions and the early-stop heuristic yield much larger improvements when added to any algorithm (see Figure \ref{fig:performance_summary}). 

Based on the previous results, we derive a few general insights into parameter value choices for the trust regions and early stop heuristic:

\begin{itemize}
    \item \textbf{Trust region growth rate ($\delta$).} A significant gain in total solution time is observed when the trust region growth rate is high, i.e., $\delta = 4$. When the size of trust regions are increased too rapidly, the computational gains from solving smaller cut generation problems are generally reduced. 
    \item \textbf{Frequency of trust region removal ($i^*$).} The total solution time is higher when $i^* = 5$ than when $i^* = 10$ or $i^* = 20$. This observation suggests that removing trust regions too frequently may reduce the computational gains of faster cut generation, especially when trust regions are small. 
    \item \textbf{Early-stop value ($\tau$).} Cuts may be weaker when $\tau$ is decreased. For example, $\tau = 1$ is the only parameter modification that leads to higher total iteration count. This is consistent with our theory, which suggests that the gain in cut strength when using trust regions comes from the computation of the extreme points of the corresponding subregion. When $\tau$ is low, an interior point of the subregion may be returned instead of an extreme point. Nonetheless, we observe through the progression of $\tau = 20$, $\tau = 5$ (baseline) and $\tau = 1$ that the improvement in cut generation speed generally outweighs the increase in iteration count from a lower $\tau$ value.
\end{itemize}


\subsection{Takeaways from numerical results}\label{subsec:computational_summary}

The three main takeaways from our numerical results are: i) trust regions improve the strength of cuts and reduce cut generation time, particularly when trust regions are small, ii) the early-stop heuristic further reduces cut generation time, whereas the dimensionality reduction can help but comes with a potential cost of reducing the strength of cuts, iii) default parameter settings seem to work well, but additional tuning of algorithm parameters can lead to significant improvements in solution time.

\section{Conclusion}

In this paper, we develop a novel class of cutting plane algorithms for solving inverse MILO problems. Our approach of using trust regions to speed up computation is simple but effective, and is well-supported by our insights on the optimality conditions of inverse MILO problems. Through extensive computational experiments, we demonstrate that our cut generation techniques and enhancements are highly effective across a large range of inverse MILO problems. 

\bibliographystyle{./StyleFiles/ormsv080}
\bibliography{InverseIntegerBib}

\newpage 

\begin{APPENDIX}{Proofs}

Several of the proofs in this section use the result that the inverse-feasible region of the sets $\mX$, $\conv(\mX)$, and $\mE(\mX)$ are the same, i.e., $\mC(\bhx,\mX) = \mC(\bhx, \conv(\mX)) = \mC(\bhx, \mE(\mX))$. This result was previously established in \cite{wang2009cutting}.


\proof{\textbf{Proof of Lemma \ref{lemma:rays_equiv}}}

The main result that must be proved is that for any given $\bx \in \mX$, $\mC(\bhx,\{\bx\}) = \mC(\bhx, \{\by(\lambda,\bhx,\bx)\}) \ \forall \lambda > 0$. Given that this statement is true, the result of $\mC(\bhx, \mE(\mX)) = \mC(\bhx,\bar\mE(\mX))$ comes trivially. The proof comes in two parts. we first prove that (i) for any given $\bx \in \mX$, $\mC(\bhx, \{\bx\}) \subseteq \mC(\bhx, \{\by(\lambda,\bhx,\bx)\}) \ \forall \lambda > 0$. We then prove that (ii) for any given $\bx \in \mX$ and $\lambda > 0$, $\mC(\bhx, \{\by(\lambda,\bhx,\bx)\}) \subseteq \mC(\bhx, \{\bx\})$. 

Proof: (i) Let $\bc \in \mC(\bhx,\{\bx\})$. By definition, $\bc^\top (\bhx - \bx) \leq 0$. Similarly, $\lambda \bc^\top (\bhx - \bx) \leq 0 \ \forall \lambda > 0$. Since $\lambda \bc^\top(\bhx - \bx) = \bc^\top (\bhx - \bhx - \lambda(\bx-\bhx)) = \bc^\top (\bhx - \by(\lambda,\bhx,\bx))$, the statement $\bc \in \mC(\bhx, \{\by(\lambda,\bhx,\bx)\}) \ \forall \lambda > 0$ must be true. (ii) For any $\bhx \in \mX$ and $\lambda >0$, let $\bc$ denote any vector such that $\lambda \bc^\top(\bhx - \bx) \leq 0$, i.e., $\bc \in \mC(\bhx, \{\by(\lambda,\bhx,\bx)\})$. Because $\lambda$ is positive, $\bc^\top(\bhx - \bx) \leq 0$. Thus, $\bc \in \mC(\bhx,\bx)$. \Halmos 
\endproof
\vspace{0.5cm}

\proof{\textbf{Proof of Theorem \ref{thm:certificate_set}}}

We prove by contradiction. ($\Leftarrow$) We assume $\mY(\bhx, \mX) = \mY(\bhx,\mG)$, but $\mG$ is not a generator set. This implies that there exists some $\bc$ such that either (i) $\bc \in \mC(\bhx,\mX)$ but $\bc \notin \mC(\bhx,\mG)$ or (ii) $\bc \in \mC(\bhx,\mG)$ but $\bc \notin \mC(\bhx,\mX)$. If (i) is true, then there exists a cost vector $\bc$ such that $\bc^\top(\bhx - \bx) \leq 0 \ \forall \bx \in \mX$ but $\bc^\top (\bhx - \bg) > 0$ for some $\bg \in \mG$. By definition, this $\bg$ cannot be in $\mX$, i.e., $\bg \notin \mX$, and thus $\bg \notin \conv(\mX)$. Furthermore, there cannot exist any $\bx \in \mX, \lambda > 0$ such that $\bhx + \lambda(\bx - \bhx) = \bg$. This implies that $\bg \notin \mY(\bhx, \mX)$, which is a contradiction. The same argument can be applied to show that (ii) leads to a contradiction.

($\Rightarrow$) We assume that $\mG$ is in fact a generator set, but $\mY(\bhx, \mX) \neq \mY(\bhx, \mG)$. This implies that there exists a point $\by^*$ such that either (i*) $\by^* \in \mY(\bhx,\mX)$ but $\by^* \notin \mY(\bhx, \mG)$ or (ii*) $\by^* \in \mY(\bhx, \mG)$  but $\by^* \notin \mY(\bhx,\mX)$. If (i*) is true, then there exists a seperating hyperplane between $\mY(\bhx, \mG)$ and $\by^*$. The normal vector of this hyperplane that is pointing in the direction of $\mY(\bhx,\mX)$ defines a cost vector $\bc$ such that all points $\bc^\top(\bhx - \by) \leq 0 \ \forall \by \in \mY(\bhx, \mG)$, which implies that $\bc^\top(\bhx-\bg) \leq 0 \ \forall \bg \in \mG$. On the other hand, $\by^*$ lies on the opposite side of the hyperplane, which implies that $\bc(\bhx - \by^*) > 0$. Finally, because $\by^* \in \mY(\bhx,\mX)$ and thus must be constructed by a convex combination of the rays $\{\by(\lambda,\bhx,\bx)\}_{\lambda \geq 0} \ \bx \in \mX$, there must exist at least one $\bx^* \in \mX$ such that $\bc(\bhx - \bx^*) > 0$. This contradicts the initial assumption that $\mC(\bhx, \mX) = \mC(\bhx, \mG)$. The same argument can be applied to show that (ii*) leads to a contradiction.
\Halmos 
\endproof
\vspace{0.5cm}

\proof{\textbf{Proof of Corollary \ref{cor:ff_cert_set}}} ($\Rightarrow$) Note that by the definition of $\mY(\bhx, \mX)$, $\mX \subseteq \mY(\bhx, \mX)$. If $\mG \subseteq \mX$ is a generator set, then $\mY(\bhx, \mG) = \mY(\bhx, \mX)$ by Theorem \ref{thm:certificate_set}. Thus, $\mX \subseteq \mY(\bhx, \mG)$. These statements together yield $\mG \subseteq \mX \subseteq \mY(\bhx,\mG)$.
($\Leftarrow$) Suppose $\mX \subseteq \mY(\bhx, \mG)$ for some $\mG \subseteq \mX$. This implies that for any $\bx \in \mX$, $\bx \in \mY(\bhx, \mG)$ and $\{\by(\lambda, \bhx, \bx)\}_{\lambda\geq 0} \in \mY(\bhx,\mG)$. Therefore, $\mY(\bhx,\mX) \subseteq \mY(\bhx, \mG)$. Since $\mG \subseteq \mX$, the reverse must also be true, i.e., $\mY(\bhx, \mG) \subseteq \mY(\bhx, \mX)$. Thus, $\mY(\bhx, \mG) = \mY(\bhx, \mX)$. Applying Theorem \ref{thm:certificate_set}, $\mG \subseteq \mX$ must be a generator set. \Halmos 
\endproof
\vspace{0.5cm}

\proof{\textbf{Proof of Theorem \ref{thm:sufficient_certificates}}}

($\Rightarrow$) We first prove the forward direction, i.e., that if $\mG \subseteq \mX$ is a generator set, then $\mB_{\epsilon}(\bhx) \cap \conv(\mX)$ must be a subset of $\conv(\mG \cup \{\bhx\})$. We prove by contradiction. Suppose that there exists a point $\by$ in $\mB_{\epsilon}(\bhx) \cap \conv(\mX)$ such that $\by \notin \conv(\mG \cup \{\bhx\})$. This implies that there exists a separating hyperplane between $\by$ and $\conv(\mG \cup \{\bhx\})$. The point $\bhx$ must lie in only of side of this hyperplane. Taking the normal vector of this hyperplane to be $\bc$, we must have that either $\bc^\top(\bhx - \bg) \leq 0 \ \forall \bg \in \mG$ and $\bc^\top (\bhx - \by) > 0$, or that $\bc^\top(\bhx - \bg) > 0 \ \forall \bg \in \mG$ and $\bc^\top (\bhx - \by) \leq 0$. Given that $\by \in \conv(\mX)$, this contracts the initial assumption that $\mC(\bhx, \mX) = \mC(\bhx, \conv(\mX)) = \mC(\bhx, \mG)$.

($\Leftarrow$) We now prove the reverse direction. In particular, we prove that if the convex hull of $\mG \subseteq \mX$ with $\bhx$ contains the set $\mB_{\epsilon}(\bhx) \cap \conv(\mX)$ for some $\epsilon>0$, then $\mG$ must be a forward-feasible generator set, i.e.,  $\mC(\bhx,\mX) = \mC(\bhx, \mG)$. We first observe that if $\mB_{\epsilon}(\bhx) \cap \conv(\mX) \subseteq \conv(\mG \cup \{\bhx\})$, then \[\mY(\bhx, \mB(\epsilon,\bhx) \cap \conv(\mX)) \subseteq \mY(\bhx, \conv(\mG \cup \{\bhx\})).\] We will now prove that both sides of this equation can be simplied to obtain $\mY(\bhx, \mX) \subseteq \mY(\bhx, \mG)$, at which point we can apply Corollary \ref{cor:ff_cert_set} to show that $\mG \subseteq \mX$ is a generator set. The major steps in this proof rely on the definition of the set $\mY(\bhx,.)$, and we refer to Theorem \ref{thm:certificate_set} for its definition. 

We first simplify the left side of the equation. Note that by the definition of $\mY(\bhx,.)$, which is a set polyhedral cone pointed at $\bhx$, $\mY(\bhx, \mB(\epsilon_1,\bhx) \cap \conv(\mX)) = \mY(\bhx, \mB(\epsilon_2,\bhx) \cap \conv(\mX))$ for any $\epsilon_1, \epsilon_2 > 0$. This implies that $\mY(\bhx, \mB(\epsilon,\bhx) \cap \conv(\mX)) = \mY(\bhx, \conv(\mX))$, since $\underset{\epsilon \rightarrow \infty}{\lim} \mB(\epsilon,\bhx) \cap \conv(\mX) = \conv(\mX)$. The set $\mY(\bhx,\conv(\mX))$ can be further simplified into $\mY(\bhx,\mX)$ by the definition of $\mY(\bhx,.)$. We now simplify the right side of the equation. Given that $\bhx \in \mY(\bhx, \conv(\mG \cup \{\bhx\}))$ by definition, the right side of the equation can be simplified to $\mY(\bhx, \conv(\mG))$. Furthermore, $\mY(\bhx, \conv(\mG)) = \mY(\bhx, \mG)$. 

These steps lead to the result that $\mY(\bhx, \mX) \subseteq \mY(\bhx, \mG)$. Since $\mX \subseteq \mY(\bhx, \mX)$, and we assumed $\mG \subseteq \mX$, $\mG \subseteq \mX \subseteq \mY(\bhx,\mG)$ and by Corollary \eqref{cor:ff_cert_set}, $\mG$ must be a forward-feasible generator set.
\Halmos 
\endproof

\end{APPENDIX}

\ECSwitch
\ECHead{Electronic Companion}

\section{\iz{Approximations using Inverse Linear Optimization}}\label{EC:approximations}

\iz{Let $\mX_{\text{LP}}$ denote a relaxed forward-feasible region where all integrality constraints present in $\mX = \{ \bx \; | \;  \bA \bx \geq \bb, \, \bx \in \mathbb{Z}^{n-q} \times \mathbb{R}^q\}$ are ignored. Similarly, let $\bc_{\text{LP}} \in \argmin_{\bc \in \mP} \{\norm{\bc - \bc^0} : \bc \in \mC(\bhx, \mX_{\text{LP}})\}$ denote the solution to an arbitrary inverse linear optimization model where $\mX_{\text{LP}}$ is considered instead of $\mX$. Since $\mX \subseteq \mX_{\text{LP}}$, the inverse-feasible region $\mC(\bhx, \mX_{\text{LP}})$ must be a subset of $\mC(\bhx, \mX)$, i.e., $\mC(\bhx, \mX_{\text{LP}}) \subseteq \mC(\bhx, \mX)$.}

\iz{The main drawback of using inverse linear optimization over $\mX_{\text{LP}}$ is that the set $\mC(\bhx, \mX_{\text{LP}})$ can be very different from the true inverse-feasible region $\mC(\bhx, \mX)$. This means that $\bc_{\text{LP}}$ can be very far from the (or in fact any) optimal solution $\bc^*$. Furthermore, there can be a large number of different regions $\mX_{\text{LP}}$ that define the same mixed integer set $\mX$ when the integrality constraints are considered. Specifically, for a fixed $\bA$ and $\bb$ generating $\mX = \{\bx \; | \; \bA \bx \geq \bb, \bx \in \mathbb{Z}^{n-q} \times \mathbb{R}^q\}$, there can exist many different $\hat \bA \neq \bA$ and $\hat \bb \neq \bb$ that can replace $\bA$ and $\bb$ and generate the same $\mX$. However, each of these $\hat \bA$'s and $\hat \bb$'s will generate different $\mX_{\text{LP}}$ sets. This means that the estimates of $\bc_{\text{LP}}$ are entirely dependent on the particular definition of $\bA$ and $\bb$, making $\bc_{\text{LP}}$, and the quality of $\bc_{\text{LP}}$ relative to $\bc^*$, unstable.}

\begin{figure}[h]
\centering
\begin{subfigure}[t]{0.24\textwidth}
    \centering
    \tikzset{every picture/.style={line width=0.75pt}} 

\begin{tikzpicture}[x=0.75pt,y=0.75pt,yscale=-1,xscale=1]
\path (0,300); 

\draw  (50,202) -- (150,202)(60,112) -- (60,212) (143,197) -- (150,202) -- (143,207) (55,119) -- (60,112) -- (65,119) (111,197) -- (111,207)(55,151) -- (65,151) ;
\draw   (118,214) node[anchor=east, scale=0.75]{1} (57,151) node[anchor=east, scale=0.75]{1} ;
\draw [color={rgb, 255:red, 245; green, 166; blue, 35 }  ,draw opacity=1 ]   (105,112) -- (125,225) ;
\draw [color={rgb, 255:red, 245; green, 166; blue, 35 }  ,draw opacity=1 ]   (40,140) -- (140,153) ;

\node[draw,circle,inner sep=1.2pt,fill,darkgray] at (111,202) {};

\node[draw,circle,inner sep=1.2pt,fill,darkgray] at (111,151) {};

\node[draw,circle,inner sep=1.2pt,fill,darkgray] at (60,151) {};

\node[draw,circle,inner sep=1.2pt,fill,darkgray] at (60,202) {};

\draw [color={rgb, 255:red, 128; green, 128; blue, 128 }  ,draw opacity=1 ]   (68,281) -- (106,281) ;
\draw [shift={(109,280)}, rotate = 538.6] [fill={rgb, 255:red, 128; green, 128; blue, 128 }  ,fill opacity=1 ][line width=0.08]  [draw opacity=0] (8.93,-4.29) -- (0,0) -- (8.93,4.29) -- cycle    ;
\draw [color={rgb, 255:red, 128; green, 128; blue, 128 }  ,draw opacity=1 ]   (68,281) -- (68,243) ;
\draw [shift={(69,240)}, rotate = 451.4] [fill={rgb, 255:red, 128; green, 128; blue, 128 }  ,fill opacity=1 ][line width=0.08]  [draw opacity=0] (8.93,-4.29) -- (0,0) -- (8.93,4.29) -- cycle    ;
\draw [color={rgb, 255:red, 128; green, 128; blue, 128 }  ,draw opacity=1 ]   (68,281) -- (96.2,270) ;
\draw [shift={(99,269)}, rotate = 518.8399999999999] [fill={rgb, 255:red, 128; green, 128; blue, 128 }  ,fill opacity=1 ][line width=0.08]  [draw opacity=0] (6.25,-3) -- (0,0) -- (6.25,3) -- cycle    ;
\draw [color={rgb, 255:red, 128; green, 128; blue, 128 }  ,draw opacity=1 ]   (68,281) -- (80,253.57) ;
\draw [shift={(82,250)}, rotate = 480.96] [fill={rgb, 255:red, 128; green, 128; blue, 128 }  ,fill opacity=1 ][line width=0.08]  [draw opacity=0] (6.25,-3) -- (0,0) -- (6.25,3) -- cycle    ;

\draw  [draw opacity=0][fill={rgb, 255:red, 80; green, 227; blue, 194 }  ,fill opacity=0.26 ]  (68,281) -- (80,254) -- (99,269) -- cycle ;




\end{tikzpicture}
    \vspace{-1.2cm}
    \caption{}
    \end{subfigure}
\begin{subfigure}[t]{0.24\textwidth}
    \centering
    \tikzset{every picture/.style={line width=0.75pt}} 

\begin{tikzpicture}[x=0.75pt,y=0.75pt,yscale=-1,xscale=1]
\path (0,300); 

\draw  (50,202) -- (150,202)(60,112) -- (60,212) (143,197) -- (150,202) -- (143,207) (55,119) -- (60,112) -- (65,119) (111,197) -- (111,207)(55,151) -- (65,151) ;
\draw   (118,214) node[anchor=east, scale=0.75]{1} (57,151) node[anchor=east, scale=0.75]{1} ;
\draw [color={rgb, 255:red, 245; green, 166; blue, 35 }  ,draw opacity=1 ]   (100,112) -- (135,225) ;
\draw [color={rgb, 255:red, 245; green, 166; blue, 35 }  ,draw opacity=1 ]   (40,124) -- (150,163) ;

\node[draw,circle,inner sep=1.2pt,fill,darkgray] at (111,202) {};

\node[draw,circle,inner sep=1.2pt,fill,darkgray] at (111,151) {};

\node[draw,circle,inner sep=1.2pt,fill,darkgray] at (60,151) {};

\node[draw,circle,inner sep=1.2pt,fill,darkgray] at (60,202) {};

\draw [color={rgb, 255:red, 128; green, 128; blue, 128 }  ,draw opacity=1 ]   (68,281) -- (106,281) ;
\draw [shift={(109,280)}, rotate = 538.6] [fill={rgb, 255:red, 128; green, 128; blue, 128 }  ,fill opacity=1 ][line width=0.08]  [draw opacity=0] (8.93,-4.29) -- (0,0) -- (8.93,4.29) -- cycle    ;
\draw [color={rgb, 255:red, 128; green, 128; blue, 128 }  ,draw opacity=1 ]   (68,281) -- (68,243) ;
\draw [shift={(69,240)}, rotate = 451.4] [fill={rgb, 255:red, 128; green, 128; blue, 128 }  ,fill opacity=1 ][line width=0.08]  [draw opacity=0] (8.93,-4.29) -- (0,0) -- (8.93,4.29) -- cycle    ;
\draw [color={rgb, 255:red, 128; green, 128; blue, 128 }  ,draw opacity=1 ]   (68,281) -- (96.2,267) ;
\draw [shift={(99,266)}, rotate = 518.8399999999999] [fill={rgb, 255:red, 128; green, 128; blue, 128 }  ,fill opacity=1 ][line width=0.08]  [draw opacity=0] (6.25,-3) -- (0,0) -- (6.25,3) -- cycle    ;
\draw [color={rgb, 255:red, 128; green, 128; blue, 128 }  ,draw opacity=1 ]   (68,281) -- (86,260.57) ;
\draw [shift={(89,257)}, rotate = 490.96] [fill={rgb, 255:red, 128; green, 128; blue, 128 }  ,fill opacity=1 ][line width=0.08]  [draw opacity=0] (6.25,-3) -- (0,0) -- (6.25,3) -- cycle    ;

\draw  [draw opacity=0][fill={rgb, 255:red, 80; green, 227; blue, 194 }  ,fill opacity=0.26 ]  (68,281) -- (96,267) -- (86,260) -- cycle ;





\end{tikzpicture}
    \vspace{-1.2cm}
    \caption{}
    \end{subfigure}
\begin{subfigure}[t]{0.24\textwidth}
    \centering
    \tikzset{every picture/.style={line width=0.75pt}} 

\begin{tikzpicture}[x=0.75pt,y=0.75pt,yscale=-1,xscale=1]
\path (0,300); 

\draw  (50,202) -- (150,202)(60,112) -- (60,212) (143,197) -- (150,202) -- (143,207) (55,119) -- (60,112) -- (65,119) (111,197) -- (111,207)(55,151) -- (65,151) ;
\draw   (118,214) node[anchor=east, scale=0.75]{1} (57,151) node[anchor=east, scale=0.75]{1} ;
\draw [color={rgb, 255:red, 245; green, 166; blue, 35 }  ,draw opacity=1 ]   (100,112) -- (135,225) ;
\draw [color={rgb, 255:red, 245; green, 166; blue, 35 }  ,draw opacity=1 ]   (40,124) -- (150,150) ;

\node[draw,circle,inner sep=1.2pt,fill,darkgray] at (111,202) {};

\node[draw,circle,inner sep=1.2pt,fill,darkgray] at (111,151) {};

\node[draw,circle,inner sep=1.2pt,fill,darkgray] at (60,151) {};

\node[draw,circle,inner sep=1.2pt,fill,darkgray] at (60,202) {};

\draw [color={rgb, 255:red, 128; green, 128; blue, 128 }  ,draw opacity=1 ]   (68,281) -- (106,281) ;
\draw [shift={(109,280)}, rotate = 538.6] [fill={rgb, 255:red, 128; green, 128; blue, 128 }  ,fill opacity=1 ][line width=0.08]  [draw opacity=0] (8.93,-4.29) -- (0,0) -- (8.93,4.29) -- cycle    ;
\draw [color={rgb, 255:red, 128; green, 128; blue, 128 }  ,draw opacity=1 ]   (68,281) -- (68,243) ;
\draw [shift={(69,240)}, rotate = 451.4] [fill={rgb, 255:red, 128; green, 128; blue, 128 }  ,fill opacity=1 ][line width=0.08]  [draw opacity=0] (8.93,-4.29) -- (0,0) -- (8.93,4.29) -- cycle    ;
\draw [color={rgb, 255:red, 128; green, 128; blue, 128 }  ,draw opacity=1 ]   (68,281) -- (96.2,270) ;
\draw [shift={(99,269)}, rotate = 518.8399999999999] [fill={rgb, 255:red, 128; green, 128; blue, 128 }  ,fill opacity=1 ][line width=0.08]  [draw opacity=0] (6.25,-3) -- (0,0) -- (6.25,3) -- cycle    ;





\end{tikzpicture}
    \vspace{-1.2cm}
    \caption{}
    \end{subfigure}
\begin{subfigure}[t]{0.24\textwidth}
    \centering
    \tikzset{every picture/.style={line width=0.75pt}} 

\begin{tikzpicture}[x=0.75pt,y=0.75pt,yscale=-1,xscale=1]
\path (0,300); 

\draw  (50,202) -- (150,202)(60,112) -- (60,212) (143,197) -- (150,202) -- (143,207) (55,119) -- (60,112) -- (65,119) (111,197) -- (111,207)(55,151) -- (65,151) ;
\draw   (118,214) node[anchor=east, scale=0.75]{1} (57,151) node[anchor=east, scale=0.75]{1} ;
\draw [color={rgb, 255:red, 245; green, 166; blue, 35 }  ,draw opacity=1 ]   (110,112) -- (135,225) ;
\draw [color={rgb, 255:red, 245; green, 166; blue, 35 }  ,draw opacity=1 ]   (40,124) -- (150,150) ;

\node[draw,circle,inner sep=1.2pt,fill,darkgray] at (111,202) {};

\node[draw,circle,inner sep=1.2pt,fill,darkgray] at (111,151) {};

\node[draw,circle,inner sep=1.2pt,fill,darkgray] at (60,151) {};

\node[draw,circle,inner sep=1.2pt,fill,darkgray] at (60,202) {};

\draw [color={rgb, 255:red, 128; green, 128; blue, 128 }  ,draw opacity=1 ]   (68,281) -- (106,281) ;
\draw [shift={(109,280)}, rotate = 538.6] [fill={rgb, 255:red, 128; green, 128; blue, 128 }  ,fill opacity=1 ][line width=0.08]  [draw opacity=0] (8.93,-4.29) -- (0,0) -- (8.93,4.29) -- cycle    ;
\draw [color={rgb, 255:red, 128; green, 128; blue, 128 }  ,draw opacity=1 ]   (68,281) -- (68,243) ;
\draw [shift={(69,240)}, rotate = 451.4] [fill={rgb, 255:red, 128; green, 128; blue, 128 }  ,fill opacity=1 ][line width=0.08]  [draw opacity=0] (8.93,-4.29) -- (0,0) -- (8.93,4.29) -- cycle    ;





\end{tikzpicture}
    \vspace{-1.2cm}
    \caption{}
    \end{subfigure}
    \caption{Plot of four different linear constraint models that define the same integer forward-feasible region when integrality constraints are considered. The forward-region in all four figures consists of points (0,0), (0,1), (1,0), and (1,1). The true inverse-feasible region for the point $\bhx = (1,1)$ is shown below each of the diagrams (as the cone between the two large perpendicular vectors), while the inverse-feasible regions of the relaxed forward-feasible region are highlighted in light blue.}\label{fig:response_example}
\end{figure}
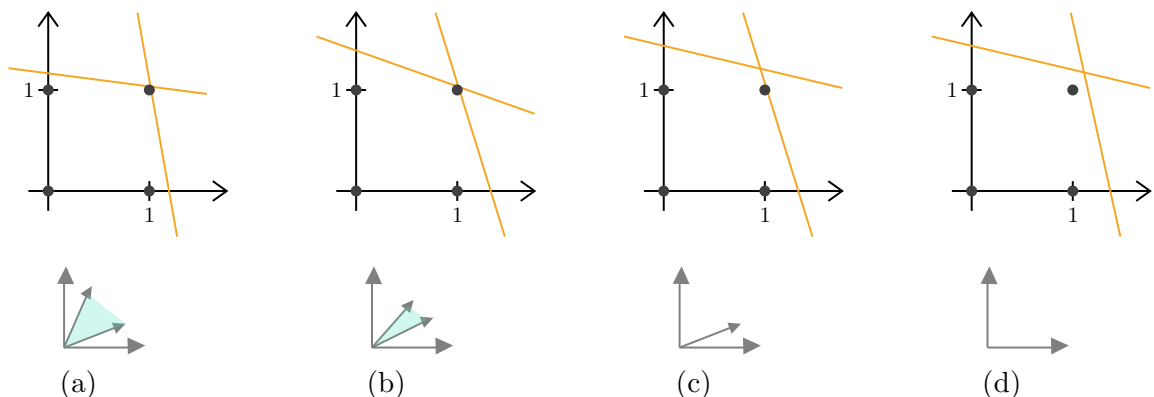

\iz{For example, consider Figure \ref{fig:response_example}, where the forward-feasible region in each of the subfigures are the same, i.e., $\mX = \{(0,0), (1,0), (0,1), (1,1)\}$. On the other hand, the linear constraints used to describe this region ($x_1, x_2 \geq 0$, and the two hyperplanes indicated by  orange lines) are different, resulting in different $\mX_{\text{LP}}$ and very different inverse-feasible regions for the solution $\bhx = (1,1)$. In particular, the inverse-feasible regions of $(\bhx, \mX_{\text{LP}})$ become increasingly smaller from Figures (1a) to (1d), which would result in solutions $\bc_{\text{LP}}$ that become increasingly worse relative to $\bc^*$ for any choice of $\bc^0$. This example, with $\mC(\bhx, \mX_{\text{LP}}) = \{\mathbf{0}\}$ in Figure (1d), also shows that there is no bound on the distance between $\bc_{\text{LP}}$ and $\bc^*$ (or similarly, between $\norm{\bc_{\text{LP}} - \bc^0}_1$ and $\norm{\bc - \bc^0}_1$). For these reasons, inverse linear optimization models may be poorly suited for solving inverse MILO models. }

\section{Example of Algorithm Performance over a Single Instance}\label{EC:example_of_instance}

\begin{figure}[h]
\centering
    \includegraphics[width = 8cm]{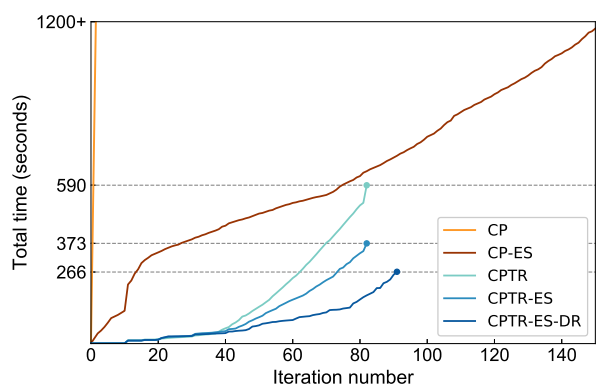}
    \caption{Total solution time over iteration count for the inverse MILO instance assign1-5-8\_t1. The instance is not solved by the CP algorithm, which manages only to compute 3 cuts within the time limit. CP-ES solves the instance in 2061 seconds using 250 iterations (not shown).}
    \label{fig:StochasticES_example}
\end{figure}

Below, we provide an example of the performance of the five algorithms over the instance assign1-5-8\_t1. In Figure \ref{fig:StochasticES_example}, the total solution time is plotted against the iteration count, and the slope of the line at any iteration reflects the cut generation speed at that particular iteration. The instance cannot be solved by CP, which has particularly high cut generation times, computing only 3 cuts by the one hour time limit. Adding the early-stop heuristic on CP (forming CP-ES) lowers cut generation times and results in a solution time of 2061 seconds with 250 iterations. On the other hand, adding trust regions on the CP (forming CPTR) lowers cut generation times more significant, \emph{and} improves the strength of cuts, resulting in the instance being solved in 590 seconds while requiring only 81 iterations. We observe through Figure \ref{fig:StochasticES_example} that cut generation times are initially very low, and gradually increase with iteration count (as a result of larger trust regions). CPTR-ES further reduces cut generation time while requiring the same number of iterations. Finally,  CPTR-ES-DR requires 9 more iterations than CPTR-ES, but the decrease in average cut generation time dominates the increase in iteration count.

\section{Table of Solved Instances}

\begin{landscape}
\small
\begin{center}
\setlength{\tabcolsep}{4.8pt}
\begin{longtable}{l r r r r r r r r r r r r}
\caption{Performance profile of all instances solved by at least one of the algorithms.} \label{TableEC:all_results} \\
& \multicolumn{2}{c}{Problem Size} & \multicolumn{2}{c}{CP } & \multicolumn{2}{c}{CP-ES} & \multicolumn{2}{c}{CPTR} & \multicolumn{2}{c}{CPTR-ES} & \multicolumn{2}{c}{CPTR-ES-DR}\\ \cmidrule(lr){2-3}\cmidrule(lr){4-5}\cmidrule(lr){6-7}\cmidrule(lr){8-9}\cmidrule(lr){10-11}\cmidrule(lr){12-13}
Instance & vars. & cons. & num.it. & time(s) &  num.it. & time(s) & num.it. & time(s) & num.it. & time(s) & num.it. & time(s)
\endfirsthead

& \multicolumn{2}{c}{Problem Size} & \multicolumn{2}{c}{CP } & \multicolumn{2}{c}{CP-ES} & \multicolumn{2}{c}{CPTR} & \multicolumn{2}{c}{CPTR-ES} & \multicolumn{2}{c}{CPTR-ES-DR}\\ \cmidrule(lr){2-3}\cmidrule(lr){4-5}\cmidrule(lr){6-7}\cmidrule(lr){8-9}\cmidrule(lr){10-11}\cmidrule(lr){12-13}
Instance & vars. & cons. & num.it. & time(s) &  num.it. & time(s) & num.it. & time(s) & num.it. & time(s) & num.it. & time(s)\\
\hline
\endhead

\hline \multicolumn{13}{|r|}{{Continued on next page}} \\ \hline
\endfoot

\hline \hline
\endlastfoot

 \toprule

50v-10\_t1 & 2013 & 233 & 0 & $>$3600 & 4964 & $>$3600 & 1567 & 1919 & 1265 & 206 & 2460 & 370 \\
50v-10\_t2 & 2013 & 233 & 0 & $>$3600 & 4567 & $>$3600 & 1002 & 992 & 990 & 121 & 1731 & 240 \\
50v-10\_t3 & 2013 & 233 & 0 & $>$3600 & 4732 & $>$3600 & 1456 & 237 & 1351 & 164 & 3046 & 402 \\
assign1-5-8\_t1 & 156 & 161 & 2 & $>$3600 & 250 & 2061 & 81 & 590 & 81 & 372 & 90 & 266 \\
assign1-5-8\_t2 & 156 & 161 & 4 & $>$3600 & 71 & 455 & 26 & 38 & 26 & 11 & 26 & 14 \\
assign1-5-8\_t3 & 156 & 161 & 3 & $>$3600 & 139 & 682 & 50 & 46 & 50 & 47 & 69 & 290 \\
bppc4-08\_t1 & 1456 & 111 & 0 & $>$3600 & 1 & 108 & 0 & $>$3600 & 1 & 108 & 1 & 108 \\
bppc4-08\_t2 & 1456 & 111 & 0 & $>$3600 & 1 & 13 & 0 & $>$3600 & 1 & 11 & 1 & 11 \\
bppc4-08\_t3 & 1456 & 111 & 0 & $>$3600 & 1 & 5 & 0 & $>$3600 & 1 & 5 & 1 & 6 \\
cod105\_t1 & 1024 & 1024 & 5 & $>$3600 & 6 & $>$3600 & 10 & $>$3600 & 17 & 1894 & 27 & 2557 \\
cod105\_t2 & 1024 & 1024 & 0 & 127 & 0 & 121 & 0 & 122 & 0 & 100 & 0 & 117 \\
cod105\_t3 & 1024 & 1024 & 6 & $>$3600 & 9 & $>$3600 & 10 & $>$3600 & 18 & 1474 & 24 & 1916 \\
csched007\_t2 & 1758 & 351 & 3 & $>$3600 & 256 & $>$3600 & 20 & $>$3600 & 147 & $>$3600 & 144 & 3475 \\
csched008\_t1 & 1536 & 351 & 17 & 1533 & 23 & 260 & 4 & 666 & 5 & 64 & 6 & 91 \\
csched008\_t2 & 1536 & 351 & 17 & 2137 & 115 & 2122 & 10 & 1489 & 40 & 2748 & 41 & 1164 \\
csched008\_t3 & 1536 & 351 & 18 & $>$3600 & 193 & 2768 & 1 & $>$3600 & 26 & 130 & 53 & 2711 \\
cvs16r128-89\_t2 & 3472 & 4633 & 0 & $>$3600 & 489 & $>$3600 & 6 & $>$3600 & 200 & 1260 & 270 & 963 \\
cvs16r128-89\_t3 & 3472 & 4633 & 0 & $>$3600 & 483 & $>$3600 & 5 & $>$3600 & 170 & 1368 & 270 & 1680 \\
drayage-100-23\_t1 & 11090 & 4630 & 660 & $>$3600 & 788 & $>$3600 & 174 & 245 & 174 & 274 & 195 & 242 \\
drayage-100-23\_t2 & 11090 & 4630 & 935 & $>$3600 & 755 & 2370 & 124 & 218 & 124 & 247 & 150 & 235 \\
drayage-100-23\_t3 & 11090 & 4630 & 492 & $>$3600 & 796 & $>$3600 & 160 & $>$3600 & 176 & $>$3600 & 178 & 3520 \\
drayage-25-23\_t1 & 11090 & 4630 & 765 & $>$3600 & 861 & $>$3600 & 162 & 277 & 159 & 270 & 173 & 217 \\
drayage-25-23\_t2 & 11090 & 4630 & 912 & 3303 & 875 & 2657 & 129 & 283 & 129 & 321 & 116 & 217 \\
drayage-25-23\_t3 & 11090 & 4630 & 461 & $>$3600 & 821 & $>$3600 & 137 & 1069 & 139 & 659 & 172 & 684 \\
eil33-2\_t1 & 4516 & 32 & 89 & 75 & 89 & 71 & 93 & 102 & 104 & 107 & 102 & 102 \\
eil33-2\_t2 & 4516 & 32 & 86 & 101 & 86 & 88 & 87 & 127 & 97 & 142 & 116 & 161 \\
eil33-2\_t3 & 4516 & 32 & 51 & 105 & 52 & 105 & 67 & 405 & 80 & 386 & 93 & 661 \\
enlight\_hard\_t1 & 200 & 100 & 0 & 0 & 0 & 0 & 0 & 0 & 0 & 0 & 0 & 0 \\
enlight\_hard\_t2 & 200 & 100 & 0 & 0 & 0 & 0 & 0 & 0 & 0 & 0 & 0 & 0 \\
enlight\_hard\_t3 & 200 & 100 & 0 & 0 & 0 & 0 & 0 & 0 & 0 & 0 & 0 & 0 \\
gen-ip002\_t1 & 41 & 24 & 63 & 704 & 60 & 229 & 44 & 600 & 44 & 90 & 53 & 105 \\
gen-ip002\_t2 & 41 & 24 & 28 & 615 & 28 & 51 & 34 & 559 & 34 & 44 & 44 & 52 \\
gen-ip002\_t3 & 41 & 24 & 65 & 893 & 61 & 309 & 52 & 688 & 54 & 144 & 70 & 189 \\
germanrr\_t1 & 10813 & 10779 & 0 & $>$3600 & 1 & 24 & 0 & $>$3600 & 1 & 21 & 1 & 23 \\
germanrr\_t2 & 10813 & 10779 & 0 & $>$3600 & 1 & 16 & 0 & $>$3600 & 1 & 22 & 1 & 22 \\
germanrr\_t3 & 10813 & 10779 & 0 & $>$3600 & 1 & 17 & 0 & $>$3600 & 1 & 22 & 1 & 22 \\
glass-sc\_t1 & 214 & 6119 & 0 & $>$3600 & 292 & $>$3600 & 161 & 961 & 161 & 140 & 188 & 151 \\
glass-sc\_t2 & 214 & 6119 & 0 & $>$3600 & 512 & 1826 & 125 & 1314 & 125 & 126 & 125 & 125 \\
glass-sc\_t3 & 214 & 6119 & 0 & $>$3600 & 640 & 2074 & 131 & 1704 & 136 & 146 & 136 & 145 \\
gmu-35-40\_t1 & 1205 & 424 & 593 & 934 & 411 & 176 & 93 & 1323 & 96 & 137 & 98 & 143 \\
gmu-35-50\_t2 & 1919 & 435 & 442 & $>$3600 & 621 & 292 & 86 & 1270 & 86 & 45 & 87 & 43 \\
gmu-35-50\_t3 & 1919 & 435 & 3 & $>$3600 & 699 & 284 & 10 & $>$3600 & 87 & 74 & 102 & 89 \\
leo1\_t1 & 6731 & 593 & 1 & 359 & 1 & 5 & 1 & 351 & 1 & 5 & 1 & 6 \\
leo1\_t2 & 6731 & 593 & 1 & 395 & 1 & 5 & 1 & 351 & 1 & 5 & 1 & 7 \\
leo1\_t3 & 6731 & 593 & 0 & 198 & 0 & 203 & 0 & 182 & 0 & 209 & 0 & 207 \\
leo2\_t1 & 11100 & 593 & 1 & 919 & 1 & 6 & 1 & 887 & 1 & 6 & 1 & 7 \\
leo2\_t2 & 11100 & 593 & 1 & 1034 & 1 & 5 & 1 & 885 & 1 & 6 & 1 & 7 \\
leo2\_t3 & 11100 & 593 & 0 & 496 & 1 & 5 & 0 & 432 & 1 & 11 & 1 & 11 \\
markshare\_4\_0\_t1 & 34 & 4 & 4 & 0 & 4 & 0 & 4 & 0 & 4 & 0 & 4 & 0 \\
markshare\_4\_0\_t2 & 34 & 4 & 4 & 0 & 4 & 0 & 4 & 0 & 4 & 0 & 4 & 0 \\
markshare\_4\_0\_t3 & 34 & 4 & 8 & 2 & 8 & 2 & 9 & 1 & 9 & 2 & 9 & 2 \\
markshare2\_t1 & 74 & 7 & 0 & $>$3600 & 8 & 7 & 0 & $>$3600 & 8 & 5 & 8 & 5 \\
markshare2\_t2 & 74 & 7 & 0 & $>$3600 & 7 & 7 & 0 & $>$3600 & 7 & 5 & 7 & 5 \\
markshare2\_t3 & 74 & 7 & 0 & $>$3600 & 8 & 8 & 0 & $>$3600 & 8 & 5 & 8 & 5 \\
mas74\_t1 & 151 & 13 & 21 & 445 & 1 & 5 & 1 & 0 & 1 & 0 & 1 & 0 \\
mas74\_t2 & 151 & 13 & 21 & 467 & 1 & 5 & 1 & 0 & 1 & 0 & 1 & 0 \\
mas74\_t3 & 151 & 13 & 21 & 458 & 1 & 5 & 1 & 0 & 1 & 0 & 1 & 0 \\
mas76\_t1 & 151 & 12 & 16 & 20 & 1 & 21 & 1 & 0 & 1 & 0 & 1 & 0 \\
mas76\_t2 & 151 & 12 & 16 & 21 & 1 & 21 & 1 & 0 & 1 & 0 & 1 & 0 \\
mas76\_t3 & 151 & 12 & 16 & 21 & 1 & 22 & 1 & 0 & 1 & 0 & 1 & 0 \\
mcsched\_t1 & 1747 & 2107 & 1 & 90 & 1 & 7 & 1 & 88 & 1 & 7 & 1 & 7 \\
mcsched\_t2 & 1747 & 2107 & 1 & 93 & 1 & 7 & 1 & 88 & 1 & 7 & 1 & 7 \\
mcsched\_t3 & 1747 & 2107 & 1 & 135 & 1 & 11 & 1 & 131 & 1 & 11 & 1 & 11 \\
mik-250-20-75-4\_t1 & 270 & 195 & 2626 & 580 & 2626 & 624 & 121 & 4 & 131 & 7 & 124 & 7 \\
mik-250-20-75-4\_t2 & 270 & 195 & 2336 & $>$3600 & 2312 & $>$3600 & 154 & 145 & 161 & 161 & 170 & 171 \\
mik-250-20-75-4\_t3 & 270 & 195 & 2480 & 451 & 2480 & 484 & 191 & 13 & 185 & 10 & 191 & 12 \\
mzzv11\_t1 & 10240 & 9499 & 616 & $>$3600 & 640 & $>$3600 & 109 & 653 & 118 & 604 & 137 & 568 \\
mzzv11\_t2 & 10240 & 9499 & 765 & $>$3600 & 732 & 3181 & 104 & 423 & 99 & 404 & 112 & 453 \\
mzzv11\_t3 & 10240 & 9499 & 723 & $>$3600 & 689 & 3318 & 132 & 698 & 156 & 815 & 152 & 527 \\
mzzv42z\_t1 & 11717 & 10460 & 764 & 2624 & 777 & 3027 & 86 & 258 & 93 & 292 & 112 & 259 \\
mzzv42z\_t2 & 11717 & 10460 & 740 & 2208 & 788 & 2326 & 100 & 186 & 103 & 257 & 111 & 186 \\
mzzv42z\_t3 & 11717 & 10460 & 953 & $>$3600 & 945 & $>$3600 & 80 & 230 & 90 & 285 & 102 & 314 \\
n5-3\_t1 & 2550 & 1062 & 8994 & 2167 & 7731 & 1889 & 1432 & 245 & 737 & 121 & 816 & 213 \\
n5-3\_t2 & 2550 & 1062 & 5645 & 1301 & 4945 & 982 & 396 & 108 & 842 & 138 & 595 & 173 \\
n5-3\_t3 & 2550 & 1062 & 5866 & 1159 & 5544 & 988 & 397 & 99 & 2355 & 461 & 495 & 163 \\
neos-2657525-crna\_t2 & 524 & 342 & 0 & $>$3600 & 19 & 116 & 0 & $>$3600 & 18 & 55 & 18 & 55 \\
neos-3083819-nubu\_t1 & 8644 & 4725 & 403 & 1265 & 319 & 1009 & 180 & 537 & 162 & 533 & 156 & 304 \\
neos-3083819-nubu\_t2 & 8644 & 4725 & 116 & 1768 & 273 & 1977 & 116 & 1094 & 218 & 1967 & 145 & 777 \\
neos-3083819-nubu\_t3 & 8644 & 4725 & 364 & 2854 & 250 & 1837 & 165 & 244 & 188 & 343 & 193 & 243 \\
neos-3381206-awhea\_t2 & 2375 & 479 & 104 & $>$3600 & 1490 & $>$3600 & 119 & $>$3600 & 379 & 3597 & 330 & 3177 \\
neos-3627168-kasai\_t1 & 1462 & 1655 & 0 & $>$3600 & 3 & 36 & 0 & $>$3600 & 3 & 36 & 4 & 11 \\
neos-3627168-kasai\_t2 & 1462 & 1655 & 0 & $>$3600 & 3 & 33 & 0 & $>$3600 & 3 & 34 & 4 & 11 \\
neos-3627168-kasai\_t3 & 1462 & 1655 & 0 & $>$3600 & 3 & 69 & 0 & $>$3600 & 3 & 64 & 3 & 34 \\
neos-4338804-snowy\_t1 & 1344 & 1701 & 0 & $>$3600 & 614 & 105 & 36 & 2 & 36 & 2 & 39 & 6 \\
neos-4338804-snowy\_t2 & 1344 & 1701 & 0 & $>$3600 & 1079 & 652 & 37 & 8 & 37 & 8 & 62 & 9 \\
neos-4338804-snowy\_t3 & 1344 & 1701 & 0 & $>$3600 & 1453 & 1391 & 95 & 50 & 95 & 51 & 143 & 60 \\
neos-4954672-berkel\_t1 & 1533 & 1848 & 0 & $>$3600 & 1909 & $>$3600 & 10 & $>$3600 & 455 & 258 & 476 & 322 \\
neos-4954672-berkel\_t2 & 1533 & 1848 & 0 & $>$3600 & 3511 & $>$3600 & 10 & $>$3600 & 424 & 205 & 462 & 215 \\
neos-4954672-berkel\_t3 & 1533 & 1848 & 0 & $>$3600 & 1423 & $>$3600 & 10 & $>$3600 & 549 & 1601 & 567 & 1547 \\
neos-860300\_t1 & 1385 & 850 & 20 & 139 & 16 & 58 & 20 & 101 & 24 & 86 & 33 & 192 \\
neos-860300\_t2 & 1385 & 850 & 19 & 122 & 20 & 66 & 31 & 165 & 35 & 144 & 40 & 162 \\
neos-860300\_t3 & 1385 & 850 & 17 & 185 & 28 & 180 & 22 & 172 & 25 & 146 & 24 & 193 \\
neos17\_t2 & 535 & 486 & 2457 & $>$3600 & 675 & 1770 & 65 & 252 & 139 & 1612 & 144 & 238 \\
neos17\_t3 & 535 & 486 & 441 & $>$3600 & 536 & 3221 & 150 & 974 & 155 & 1958 & 152 & 556 \\
neos5\_t1 & 63 & 63 & 6 & $>$3600 & 69 & 688 & 38 & 0 & 38 & 0 & 38 & 0 \\
neos5\_t2 & 63 & 63 & 71 & 1318 & 67 & 470 & 36 & 0 & 36 & 0 & 36 & 0 \\
neos5\_t3 & 63 & 63 & 84 & 3343 & 81 & 2760 & 38 & 0 & 38 & 0 & 38 & 0 \\
qap10\_t1 & 4150 & 1820 & 19 & 514 & 35 & 568 & 24 & 641 & 34 & 579 & 34 & 585 \\
qap10\_t2 & 4150 & 1820 & 23 & 461 & 32 & 428 & 30 & 571 & 39 & 543 & 39 & 542 \\
qap10\_t3 & 4150 & 1820 & 27 & 451 & 27 & 264 & 37 & 854 & 51 & 1022 & 51 & 1063 \\
ran14x18-disj-8\_t1 & 504 & 447 & 0 & $>$3600 & 994 & 726 & 625 & 180 & 662 & 96 & 731 & 96 \\
ran14x18-disj-8\_t2 & 504 & 447 & 0 & $>$3600 & 1777 & 1289 & 569 & 88 & 569 & 95 & 685 & 103 \\
ran14x18-disj-8\_t3 & 504 & 447 & 0 & $>$3600 & 1250 & 788 & 563 & 110 & 522 & 60 & 667 & 71 \\
rocI-4-11\_t1 & 6839 & 10883 & 1 & 268 & 1 & 296 & 1 & 321 & 1 & 369 & 1 & 368 \\
rocI-4-11\_t2 & 6839 & 10883 & 2 & 368 & 2 & 318 & 2 & 681 & 2 & 668 & 2 & 645 \\
rocI-4-11\_t3 & 6839 & 10883 & 3 & 226 & 5 & 83 & 3 & 146 & 4 & 56 & 4 & 63 \\
rococoB10-011000\_t1 & 4456 & 1667 & 0 & $>$3600 & 1 & 9 & 0 & $>$3600 & 1 & 9 & 1 & 9 \\
rococoB10-011000\_t2 & 4456 & 1667 & 0 & $>$3600 & 1 & 9 & 0 & $>$3600 & 1 & 9 & 1 & 9 \\
rococoB10-011000\_t3 & 4456 & 1667 & 0 & $>$3600 & 1 & 11 & 0 & $>$3600 & 1 & 10 & 1 & 10 \\
rococoC10-001000\_t1 & 3117 & 1293 & 1 & 261 & 1 & 10 & 1 & 335 & 1 & 9 & 1 & 10 \\
rococoC10-001000\_t2 & 3117 & 1293 & 1 & 269 & 1 & 10 & 1 & 335 & 1 & 9 & 1 & 10 \\
rococoC10-001000\_t3 & 3117 & 1293 & 1 & 270 & 1 & 10 & 1 & 335 & 1 & 9 & 1 & 10 \\
roi2alpha3n4\_t1 & 6816 & 1251 & 20 & $>$3600 & 154 & $>$3600 & 149 & $>$3600 & 146 & 1599 & 142 & 1455 \\
roi2alpha3n4\_t2 & 6816 & 1251 & 21 & $>$3600 & 197 & $>$3600 & 130 & $>$3600 & 125 & 524 & 144 & 3294 \\
roi2alpha3n4\_t3 & 6816 & 1251 & 20 & $>$3600 & 155 & $>$3600 & 117 & $>$3600 & 146 & 2830 & 148 & 3537 \\
roll3000\_t1 & 1166 & 2295 & 1 & 9 & 1 & 10 & 1 & 13 & 2 & 23 & 2 & 20 \\
roll3000\_t2 & 1166 & 2295 & 10 & 168 & 5 & 81 & 5 & 83 & 5 & 72 & 5 & 80 \\
roll3000\_t3 & 1166 & 2295 & 1 & 7 & 1 & 7 & 1 & 0 & 1 & 5 & 1 & 5 \\
seymour\_t1 & 1372 & 4944 & 0 & $>$3600 & 443 & $>$3600 & 10 & $>$3600 & 1136 & 449 & 1426 & 431 \\
seymour\_t2 & 1372 & 4944 & 0 & $>$3600 & 474 & $>$3600 & 10 & $>$3600 & 959 & 272 & 1036 & 253 \\
seymour\_t3 & 1372 & 4944 & 0 & $>$3600 & 450 & $>$3600 & 30 & $>$3600 & 1194 & 614 & 1532 & 643 \\
seymour1\_t1 & 1372 & 4944 & 47 & $>$3600 & 331 & $>$3600 & 1679 & 800 & 1543 & 521 & 1663 & 454 \\
seymour1\_t2 & 1372 & 4944 & 48 & $>$3600 & 396 & $>$3600 & 956 & 495 & 970 & 263 & 1000 & 238 \\
seymour1\_t3 & 1372 & 4944 & 49 & $>$3600 & 410 & $>$3600 & 1669 & 920 & 1669 & 496 & 1875 & 498 \\
sp150x300d\_t1 & 600 & 450 & 4424 & $>$3600 & 1414 & 133 & 248 & 35 & 248 & 36 & 205 & 24 \\
sp150x300d\_t2 & 600 & 450 & 2035 & 181 & 2035 & 181 & 222 & 26 & 222 & 29 & 228 & 19 \\
sp150x300d\_t3 & 600 & 450 & 1653 & 193 & 1606 & 154 & 199 & 18 & 199 & 19 & 212 & 27 \\
splice1k1\_t1 & 3253 & 6505 & 0 & $>$3600 & 1 & 48 & 0 & $>$3600 & 1 & 40 & 1 & 50 \\
splice1k1\_t2 & 3253 & 6505 & 0 & $>$3600 & 1 & 47 & 0 & $>$3600 & 1 & 49 & 1 & 41 \\
splice1k1\_t3 & 3253 & 6505 & 0 & $>$3600 & 1 & 48 & 0 & $>$3600 & 1 & 40 & 1 & 50 \\
supportcase26\_t1 & 436 & 870 & 38 & 729 & 36 & 37 & 16 & 3 & 16 & 3 & 16 & 3 \\
supportcase26\_t2 & 436 & 870 & 13 & $>$3600 & 70 & 2018 & 29 & 1537 & 30 & 686 & 34 & 950 \\
supportcase26\_t3 & 436 & 870 & 9 & 487 & 14 & 8 & 3 & 0 & 3 & 0 & 3 & 0 \\
wachplan\_t1 & 3361 & 1553 & 0 & 1787 & 0 & 1845 & 0 & 1123 & 0 & 1087 & 0 & 1087 \\
wachplan\_t2 & 3361 & 1553 & 0 & 1804 & 0 & 1823 & 0 & 1122 & 0 & 1131 & 0 & 1135 \\
wachplan\_t3 & 3361 & 1553 & 0 & 1799 & 0 & 1823 & 0 & 1123 & 0 & 1089 & 0 & 1133\\
\end{longtable}
\end{center}
\end{landscape}


\end{document}